\documentstyle[11pt]{article}

\newcommand{\qed}{\hfill$\Box$}

\newcommand{\D}{{\it Des}}
\newcommand{\des}{{\it des}}

\newcommand{\del}{{\it del}}

\newcommand{\rmaj}{{\it rmaj}}

\newtheorem{thm}{Theorem}[section]
\newtheorem{pro}[thm]{Proposition}
\newtheorem{lem}[thm]{Lemma}

\newtheorem{cor}[thm]{Corollary}

\newtheorem{exa}[thm]{Example}
\newtheorem{df}[thm]{Definition}
\newtheorem{rem}[thm]{Remark}

\def\sg{\sigma}

\def\kap{\kappa}

\begin{document}

\title{$q$ Statistics on $S_n$ and Pattern Avoidance}

\author{Amitai Regev\thanks{Partially supported by Minerva
Grant No. 8441 and by EC's IHRP Programme, within the Research Training Network ``Algebraic Combinatorics
in Europe'', grant HPRN-CT-2001-00272}\\
Department of Mathematics\\
The Weizmann Institute of Science\\
Rehovot 76100, Israel\\
{\em regev@wisdom.weizmann.ac.il}\\
\and
Yuval Roichman\thanks{Partially supported by EC's IHRP Programme, within the Research Training Network ``Algebraic Combinatorics
in Europe'', grant HPRN-CT-2001-00272}\\
Department of Mathematics\\
Bar Ilan University\\
Ramat Gan 52900, Israel\\
{\em yuvalr@math.biu.ac.il}}
\date{May 28, 2003}

\maketitle

\bibliographystyle{acm}

\begin{abstract}
Natural $q$ analogues of classical statistics 
on the symmetric groups $S_n$ are introduced; parameters like: the $q$-length, 
the $q$-inversion number, the $q$-descent number and the $q$-major index.
MacMahon's theorem about the 
equi-distribution of 
the inversion number and the reverse major index
is generalized to all 
positive integers $q$.
It is also shown that the $q$-inversion number and the $q$-reverse major index are equi-distributed over subsets of permutations
avoiding certain patterns. 
Natural $q$ analogues of the Bell and the Stirling numbers
are related to these $q$ statistics
 - through the counting of the above pattern-avoiding permutations. 
\end{abstract}

\section{Introduction}

MacMahon's celebrated theorem about the equi-distribution of the {\it length}
(or the {\it inversion-number})
and the {\it major index} statistics on the symmetric group 
$S_n$ -- has received
far--reaching refinements and generalizations through the last three decades. 
For a brief review on these refinements -- see~\cite{RR}.
In [15] we extended the various classical $S_n$ statistics, in a natural way, 
to the alternating group $A_{n+1}$. This was done via the canonical
presentations of the elements of these groups, and by a certain covering map
$f:A_{n+1}\to S_n$.

Further refinements of MacMahon's theorem
were obtained in [15] by the introduction of the 
{\it `delent'} statistics on these groups. Then
these equi-distribution theorems for $S_n$ were `lifted' back,
via $f:A_{n+1}\to S_n$, thus yielding equi-distribution theorems 
for $A_{n+1}$.

This paper continues \cite{RR} and might be considered as its {\it q-analogue}.
We introduce the {\it q-analogues} of the classical statistics on the symmetric
groups: the {\it q-length}, the {\it q-inversion number}, 
the {\it q-descent number},
the {\it q-major index} and the {\it q-reverse-major index} of a 
permutation. The {\it q-delent} statistics are also introduced. We
then extend  classical  properties to these $q$ -analogues. For example,
it is proved that the $q$--length equals the $q$--inversion number of a 
permutation; furthermore, it is proved that the $q$-inversion number and the $q$-reverse major index are equidistributed on $S_{n+q-1}$. See below.

It is realized that the above map $f:A_{n+1}\to S_n$
is the restriction to $A_{n+1}$ 
of a covering map $f_2:S_{n+1}\to S_n$. More generally, we have similar
covering maps 
$f_q:S_{n+q-1}\to S_n$ for all positive
integers $q$. These maps are defined via the 
canonical presentations of the elements in $S_{n+q-1}$. It is proved 
that the map $f_q$ sends the $q$--statistics on $S_{n+q-1}$ to
the corresponding
classical statistics on $S_n$, see Proposition~\ref{cover1}
below. For example, if $\pi\in S_{n+q-1}$, it is proved there that
the $q$--inversion number of $\pi$ equals the inversion number
of $f_q(\pi)$. 

\bigskip

Dashed patterns in permutations were introduced by Babson and  Steingrimsson~\cite{BS}. For
example, a permutation $\sg$ contains the pattern $(1-32)$ if 
$\sg=[\ldots,a,\ldots,c,b,\ldots ]$ for some $a<b<c$; if no such
$a,b,c$ exist then $\sg$ is said to avoid  $(1-32)$. 
Connections between the number of
permutations avoiding $(1-32)$ -- and various combinatorial objects, like 
the Bell and the Stirling numbers, as well as the number of
left-to-right-minima in permutations were proved by Claesson~\cite{Cl} . 
Via the various $q$--statistics we obtain $q$-analogues for these connections 
and results. 

For a permutation $\pi\in S_{n+q-1}$ it is proved that 
the {\it q--descent} and the {\it q--delent} numbers of $\pi$ are equal 
exactly when $\pi$  avoids a certain collection of dashed patterns, 
and that the number of these permutations is $(q-1)!\sum_k q^k S(n,k)$,
where $S(n,k)$ are the Stirling numbers of the second kind,
see Proposition~\ref{main08}. 
Also, the number of permutations in $S_{n+q-1}$ for which the {\it q--delent}
number equals $k-1$ is $(q-1)!q^kc(n,k)$, where $c(n,k)$ are the Stirling numbers
of the first kind; see Proposition~\ref{main09}.

\bigskip

Equi-distribution of $q$-statistics is studied in Section~\ref{equid}.
A {\it q-analogue} of MacMahon's  
classical equi-distribution theorem is given,
see Theorems~\ref{main05} 
below.
Multivariate refinements of MacMahon's theorem, due to 
Foata-Sch\"utzenberger and others \cite{Ros, FS, RR},
have also corresponding $q$-analogues.
These analogues are described in Section~\ref{q-MM}, see 
Theorem~\ref{fs-q2} and its consequences.

An intensive study of equi-distribution 
over subsets of permutations avoiding patterns 
has been carried out recently,
cf.~\cite{E, EP, Rei, AR}.
In Section~\ref{equid-avoid} it is shown that certain $q$-statistics are
equi-distributed on the above 
subsets of dashed--patterns--avoiding permutations.
See Theorem~\ref{main06} and~\ref{q6} below.

\section{The main results}\label{mr}

Throughout the paper $q$ is a positive integer.
Recall 
the unique canonical
presentation of a permutation in $S_n$
as a product of shortest coset representatives along the principal flag,
see Subsection~\ref{Sn} below.
The {\it $q$-length} of a permutation $\pi\in S_n$, 
$\ell_q(\pi)$, is the number of 
Coxeter generators
in the canonical presentation of $\pi$, where  
the generators $s_1,\dots,s_{q-1}$ are not counted.
~Define the {\it $q$-inversion number} $inv_q(\pi )$ as
$$
inv_q(\pi):=0\quad\mbox{if $n\le q$; and if $q<n$},
$$

$$
inv_q(\pi):=\sum\limits_{i=q+1}^n m_q(i),
$$
where 
$$
m_q(i):=\min\{i-q,\#\{j<i|\pi(j)>\pi(i)\} \}. 
$$
%
Thus $\ell_1(\pi)=\ell(\pi)$   and $inv_1(\pi)=inv(\pi)$.

\medskip

As in the (classical) case $q=1$, we have
\begin{pro}\label{main01} [See Proposition~\ref{cover01}].
For every $\sg \in S_n$
$$
\ell_q(\sg)=inv_q(\sg).
$$
\end{pro}


\medskip

\begin{pro}\label{main001}[ See Proposition~\ref{ex3}].
For every 
$\pi\in A_{n}$,
$\ell_2(\pi)$ is the length with respect to the 
set of generators 
$\{a_1,\dots,a_{n-2}\}\subset A_n$,
where $a_i:=s_1s_{i+1}$.
\end{pro}


\medskip

Define the {\it $q$-delent number}, $del_q(\pi)$, to be
the number of times $s_q$ appears in the canonical presentation of $\pi$.

For $\,0\le k\le n-1$
define the {\it $k$-th almost left to right minima}
in a permutation $\pi\in S_n$ (denoted a$^k$.l.t.r.min)
as the set of indices
$$
Del_{k+1}(\pi):=
\{ i \mid k+2 \le i \le n, \#\{j<i \mid \pi (j) < \pi (i) \} \le k \}.
$$
Thus $Del_q(\pi)$ is the set of a$^{q-1}$.l.t.r.min in $\pi$.
%
See Example~\ref{exa} below.

\medskip

\begin{pro}\label{main02} [See Proposition~\ref{altr2}].
The number of occurences of $s_{k+1}$ 
in the canonical presentation of $\pi\in S_n$, $del_{k+1}(w)$, 
equals the number of a$^k$.l.t.r.min in  $\pi$. 
\end{pro}

\smallskip

\noindent
The second delent statistics $del_2$
on even permutations in $A_{n+1}$ 
and the first delent statistics $del_1$ on $S_n$
have analogous interpretations. See, for example, Proposition~\ref{ex3}.

\bigskip

The {\it $q$-descent set} of $\pi\in S_{n+q-1}$ is defined as
$$
Des_q(\pi):=\{i\ |\ i \hbox{ is a $q$-descent in } \pi\},   
$$
and 
the {\it $q$-descent number} is defined as
$$
des_q(\pi):=\#Des_q(\pi). 
$$
For $\pi\in S_{n+q-1}$ define the {\it $q$-major index}
$$
maj_{q}(\pi):=\sum\limits_{i\in Des_q(\pi)} i
$$
and the {\it $q$-reverse major index}
$$
rmaj_{q,m}(\pi):=\sum\limits_{i\in Des_q(\pi)} (m-i),
$$
where $m=n+q-1$.

\medskip

Thus $Des_1$ is the standard descent set of a permutation in $S_n$. 
The definition of the $q$-descent set is justified by the following 
phenomena: 
\begin{itemize}
\item[(1)]   $Des_2$ is the descent set on 
the alternating group $A_{n}$ with respect to the distinguished set of generators $\{a_1,\dots,a_{n-2}\}$, where $a_i:=s_1s_{i+1}$,  
see Proposition~\ref{ex3}. 
\item[(2)] The $q$-descent set, $Des_q$, is strongly related with pattern avoiding
permutations, see Proposition~\ref{q-avoid}.  
\item[(3)] $Des_q$ is  involved in the definition of the
$q$-(reverse) major index, and thus in the $q$-analogue of MacMahon's equi-distribution 
theorem (Theorem~\ref{qmac}).
\end{itemize}

\bigskip

Given $q$, denote by
\[
Pat(q)=\{(\sg_1-\sg_2-\cdots-\sg_q-(q+2),(q+1)) \mid \sg\in S_q\}
\]
the set with these $q!$ dashed patterns.\\
For example, $Pat(1)=\{(1-32)\}$
$Pat(2)=\{(1-2-43),\;(2-1-43)\}$.

Denote by $Avoid_q(n+q-1)$ the set of permutations in $S_{n+q-1}$
avoiding all the $q!$ patterns in $Pat(q)$.
%

\begin{pro}\label{main04} [See Proposition~\ref{q-avoid}].
A permutation $\pi\in S_{n+q-1}$ avoids $Pat(q)$ exactly
when $Del_q(\pi)-1=Des_q(\pi)$:
$$
Avoid_q(n+q-1)=\{\pi\in S_{n+q-1}|\ Del_q(\pi)-1=Des_q(\pi)\}.
$$
\end{pro}

\bigskip

The following is a $q$-analogue of MacMahon's equi-distribution theorem.

\begin{thm}\label{main05} [See theorem~\ref{qmac}].
$$
\begin{array}{ll}
\sum\limits_{\pi\in S_{n+q-1}} t^{rmaj_{q,n+q-1}(\pi)}
=
\sum\limits_{\pi\in S_{n+q-1}} t^{inv_q(\pi)}=~~~~~~~~~\\
~~~~~~~~~~~~~~~~~~~~~~=q!(1+tq)(1+t+t^2q)\cdots 
(1+t+\ldots +t^{n-2}+t^{n-1}q).
\end{array}
$$
\end{thm}

Far reaching multivariate refinements of MacMahon's theorem, which
imply equi-distribution on subsets
of permutations,  were given by 
Foata and Sch\"utenberger and others, cf.~\cite{Ros, FS, GG, RR}.
In Subsection~\ref{q-MM} we describe some $q$-analogues of these
refinements, see Theorem~\ref{fs-q} and Corollary~\ref{q-rosel}
below.

\medskip

The above $q$-statistics are equi-distributed on permutations avoiding $Pat(q)$.

\begin{thm}\label{main06} [See Corollary~\ref{qs5}].
$$
\sum\limits_{\pi ^{-1}\in Avoid_q(n+q-1)} t_1^{rmaj_{q,n+q-1}(\pi)}t_2^{des_q(\pi)}
=
\sum\limits_{\pi ^{-1}\in Avoid_q(n+q-1)} t_1^{inv_q(\pi)}t_2^{des_q(\pi)}
$$
\end{thm}

\smallskip

\noindent
For example, for $q=1$ 
$$
\sum\limits_{\pi ^{-1}\in Avoid(1-32)} t_1^{rmaj_{n}(\pi)}t_2^{des(\pi)}
=
\sum\limits_{\pi ^{-1}\in Avoid(1-32)} t_1^{inv(\pi)}t_2^{des(\pi)}.
$$
For $q=2$
$$
\sum\limits_{\pi ^{-1}\in Avoid(1-2-43, 2-1-43)} t_1^{rmaj_{2,n+1}(\pi)}t_2^{des_2(\pi)}
=
\sum\limits_{\pi ^{-1}\in Avoid(1-2-43, 2-1-43)} t_1^{inv_2(\pi)}t_2^{des_2(\pi)}.
$$


\bigskip

Bell  and Stirling numbers (of both kinds)
appear naturally in the enumeration of permutations with respect to
their $q$-statistics. 

Let $c(n,k)$ be the $k$-th Stirling number of the first kind
and $S(n,k)$ be the $k$-th Stirling number of the second kind.
Let the $n$-th $q$-Bell number be
$b_q(n):=\sum_k q^k S(n,k)$.
Let $B_q(x):=\sum_{n=0}^{\infty}b_q(n)\frac{x^n}{n!}$ denote the 
exponential generating function  of $\{b_q(n)\}$. Then
\[
B_q(x)=\exp (q e^x-q).
\]
The classical formula
$b_1(n)=\frac{1}{e}\sum_{r=0}^{\infty}\frac{r^n}{r!}$   ~\cite{D}
(see also~\cite[(1.6.10)]{W})
generalizes as follows: 
\[
b_q(n)=\frac{1}{e^q}\sum_{r=0}^{\infty}\frac{q^rr^n}{r!},
\]
see Remark~\ref{remarkable}.

\medskip

\noindent We have 

\begin{pro}\label{main07} [See Proposition~\ref{qpro}].
\[
\begin{array}{ll}
\#\{\sg\in S_{n+q-1}\mid Del_q(\sg)-1=Des_q(\sg)\quad\mbox{and}\quad
del_q(\sg)=k-1\} =~~~~~~~~~~~~~~~~~~\\
~~~~~~~~~~~~~~~~~~~~~~~~~~~~~~~~~~~~~~~~~~~~~~~~~~~~~~~~~~~~=(q-1)!q^kS(n,k).
\end{array}
\]
\end{pro}

\begin{cor}\label{main08}  [See Propositions~~\ref{q-avoid} 
and~\ref{nu1}].
\[
(q-1)!b_q(n)=
\#\{\pi\in S_{n+q-1}|\ Del_q(\pi)-1=Des_q(\pi)\}=Avoid_q(n+q-1).
\]
\end{cor}

\begin{pro}\label{main09} [See Propositions~\ref{qc3}].
$$
\#\{\pi\in S_{n+q-1}|\ del_q(\pi)=k-1\}=c_q(n,k),
$$
where $c_q(n,k)=q^k\ (q-1)!\ c(n,k)$.
\end{pro}


\section{Preliminaries}\label{pre} 

\subsection{The $S_n$ Canonical Presentation}\label{Sn}

A basic tool, both in~\cite{RR} and in this paper, is the canonical presentation
of a permutation, which we now describe.

\medskip

Recall that the transpositions
$s_i=(i,i+1)$, $1\le i<n-1$, are the Coxeter generators of the symmetric
group $S_n$. For each $1\le j\le n-1$ define 
\begin{eqnarray}\label{eqn1}
R^S_j=\{ 1,\,s_j,\,s_js_{j-1},\ldots , \,s_js_{j-1}\cdots s_1\}
\end{eqnarray}
and note that $R^S_1,\ldots ,R^S_{n-1} \subseteq S_n$.\\
The following is a classical theorem; see
for example~\cite[pp. 61-62]{Go}. See also~\cite[Theorem 3.1]{RR}.
\begin{thm}\label{RR-3.1} 
Let $w\in S_n$, then there exist unique elements $w_j\in R^S_j$,
 $1\le j\le n-1$, such that $w=w_1\cdots w_{n-1}$. 
Thus, the presentation $w=w_1\cdots w_{n-1}$ is unique; it 
is called the {\bf canonical presentation} of $w$. 
\end{thm}

Note that $R^S_j$ is the complete list of representatives
of minimal length of right cosets of $S_{j}$ in $S_{j+1}$.
Thus, the canonical presentation of $w\in S_n$ is the unique presentation
of $w$ as a product of shortest coset representatives along the principal flag
$$
\{e\}= S_1< S_2<\cdots <S_n.
$$

\bigskip

We remark that a similar canonical presentation for the alternating groups
$A_n$ -- is given in~\cite{RR}, see Section~\ref{An} below.

\medskip

The descent set $Des(\pi)$ of a permutation $\pi\in S_n$ is a classical
notion. In~\cite{RR} the {\it `delent'} statistic was introduced: 
$Del(\pi)$ is the set of indices $i$ which are left-to-right-minima 
of $\pi$, and $del(\pi)=\# Del(\pi)$. By Proposition 7.2 of~\cite{RR},
$del(\pi)$ equals the number of times that $s_1=(1,2)$ appears in the 
 canonical presentation of $\pi$.

Theorem 9.1 is the main theorem of~\cite{RR} and we now state its 
part about $S_n$ (it also has a similar part about $A_n$).

\begin{thm}\label{9.1}
For every subsets $D_1\subseteq [n-1]$ and $D_2\subseteq [n-1]$
$$
\sum\limits_{\{\pi\in S_n|\ Des_S(\pi^{-1})\subseteq D_1,\ 
Del_S(\pi^{-1})\subseteq D_2\}} q^{\rmaj_{S_n}(\pi)}
=
$$
$$
=\sum\limits_{\{\pi\in S_n|\ Des_S(\pi^{-1})\subseteq D_1,\
Del_S(\pi^{-1})\subseteq D_2\}} q^{\ell_S(\pi)}.
$$
%
\end{thm}

In the following case, a simple explicit generating function is given.

\begin{thm}\label{6.1}
\cite[Theorem 6.1]{RR}
$$
\sum_{\sg\in S_n} q^{\ell_S(\sg)}t^{del_S(\sg)}=
\sum_{\sg\in S_n} q^{\rmaj_{S_n}(\sg)}t^{del_S(\sg)}= 
$$
$$
=(1+qt)(1+q+q^2t)\cdots (1+q+\ldots  +q^{n-1}t). 
$$
\end{thm}

\subsection{The Alternating Group}\label{An}

The alternating group serves as a motivating example.
Here are some results from~\cite{RR}, which are applied in Sections~\ref{ex}
and~\ref{appendixI}
and in the formulation and proof of Proposition~\ref{f2}.
The reader who is not interested in this motivating example may
skip this subsection. 

Let 
$$
a_i:=s_1s_{i+1}\qquad (1\le i\le n-1).
$$
The set 
$$
A:=\{a_i\ | \ 1\le i\le n-1\}
$$
generates the alternating group on $n$ letters $A_{n+1}$.
This generating set and its following properties appear in \cite{Mi}.

\begin{pro}\label{RR-m0}\cite[Proposition 2.5]{Mi}
The defining relations of $A$ are
$$
(a_ia_j)^2=1 \qquad (|i-j|>1) ;
$$
$$
(a_ia_{i+1})^3=1 \qquad (1\le i< n-1) ;
$$
$$
a_1^3=1\qquad\mbox{and}\qquad a_i^2=1 \qquad (1<i\le n-1) .
$$
\end{pro}

For each $1\le j\le n-1$ define 
\begin{eqnarray}\label{eqn2}
R^A_j=\{ 1,\;a_j,\;a_ja_{j-1},\;\ldots ,\; a_j\cdots a_2,\;a_j
\cdots a_2a_1,\;a_j\cdots a_2a^{-1}_1\}
\end{eqnarray}
and note that $R^A_1,\ldots ,R^A_{n-1} \subseteq A_{n+1}$.

\begin{thm}\label{RR-3.4}
Let $v\in A_{n+1}$, then there exist unique elements $v_j\in R^A_j$,
 $1\le j\le n-1$, such that $v=v_1\cdots v_{n-1}$, and this presentation 
is unique.

\end{thm}

This presentation is called the $A$ canonical presentation of $v$.

For $\sg\in A_{n+1}$ let $\ell_A(\sg)$ be the length of the $A$ canonical
presentation of $\sg$.
Let 
$$
Des_A(\sg):=\{i\ |\ \ell_A(\sg)\le \ell_A(\sg a_i)\}
$$
and $des_A(\sg):=\# Des_A(\sg)$, ~define ~$maj_{A}(\sg):=
\sum\limits_{i\in Des_a(\sg)} i$, \\
and $rmaj_{A_{n+1}}(\sg):=
\sum\limits_{i\in Des_a(\sg)} (n-i)$. Let $del_A(\sg)$ the number of 
appearances of $a_1^{\pm 1}$ in its $A$ canonical presentation.
It is proved in~\cite{RR} that this number equals the number
of almost-left-to-right-minima in $\sg$.

\bigskip

Theorems~\ref{RR-3.1} and~\ref{RR-3.4} allow us to introduce in \cite{RR}
the following covering map :
\begin{df}\label{RR-df2} 
Define  $f:A_{n+1}\to S_n$ as follows. 
\[
f(a_1)=f(a^{-1}_1)=s_1\qquad\mbox{and}\qquad f(a_i)=s_i,
\qquad 2\le i\le n-1.
\]
Now extend $f:R^A_j\to R^S_j$ via
\[
f(a_ja_{j-1}\cdots a_{\ell})=s_js_{j-1}\cdots s_{\ell},\qquad
f(a_j\cdots a_1)=f(a_j\cdots a^{-1}_1)=s_j\cdots s_1.
\]
Finally, let $v\in A_{n+1},\quad v=v_1\cdots v_{n-1}$ its unique
$A$ canonical presentation, then 
\[
f(v)=f(v_1)\cdots f(v_{n-1})\,
\]
which is clearly the $S$ canonical presentation of $f(v)$.
\end{df}

\begin{pro}\label{RR-5.34}\cite[Propositions 5.3-5.4]{RR}
For every $\pi\in A_{n+1}$, 
$$
\ell_A(\pi)=\ell_S(f(\pi)),\ \ \ \ \ \ \  \D_A(\pi)=Des_S(f(\pi)),
 \ \ \ \ \ \ \
Del_A(\pi)=Del_S(f(\pi))
$$
Thus
$\des_A(\pi)=\des_S(f(\pi))$,
$maj_A(\pi)=maj_S(f(\pi))$,
$\rmaj_{A_{n+1}}(\pi)=rmaj_{S_n}(f(\pi))$ and
$\del_A(\pi)=\del_S(f(\pi))$.
\end{pro}

\section{Basic Concepts I}\label{invlg} 

Let $\pi\in S_n$. Recall that its length  $\ell (\pi)$ 
equals the number of the
Coxeter generators $s_1,\dots,s_{n-1}$
in its canonical presentation. It is well known that $\ell (\pi)$ also equals 
$inv (\pi)$, the number of inversions of $\pi$. Also, it is easily seen that 
$inv (\pi)$ can be written as 
\[
inv(\pi)=\sum\limits_{i=2}^n m(i),
\]
where 
\[
m(i)=\min\{i-1,\#\{j<i|\pi(j)>\pi(i)\} \}. 
\]
Thus, the following definition is a natural $q$-analogue of these 
two classical statistics. 

\begin{df}\label{qinvlength}
Let $\pi\in S_n$.
\begin{enumerate}
\item ($\ell_q$)
~Let $q< n$ and define the {\it $q$-length} 
$\ell_q(\pi)$ as follows:\\
$\ell_q(\pi):=$ the number of 
Coxeter generators
in the canonical presentation of $\pi$, where  
$s_1,\dots,s_{q-1}$ are not counted
(thus, for example, $\ell_2(s_1)=0$ and   $\ell_2(s_1s_2s_1s_3s_2s_1)=3$).
\item ($inv_q$)
~Define the {\it $q$-inversion number} $inv_q(\pi )$ as
$$
inv_q(\pi):=0\quad\mbox{if $n\le q$; and if $q<n$},
$$

$$
inv_q(\pi):=\sum\limits_{i=q+1}^n m_q(i),
$$
where 
$$
m_q(i):=\min\{i-q,\#\{j<i|\pi(j)>\pi(i)\} \}. 
$$
\end{enumerate}
\end{df}
Thus $\ell_1(\pi)=\ell(\pi)$   and $inv_1(\pi)=inv(\pi)$.

\medskip

As in the (classical) case $q=1$, we have
\begin{pro}\label{cover01}
For every $\sg \in S_n$
$$
\ell_q(\sg)=inv_q(\sg).
$$
\end{pro}

\noindent
{\bf Proof.}  We may assume that $q<n$.
Let $\sg=w_1\cdots w_{n-1}$ with $w_j\in R_j$ be the canonical
presentation of $\sg$, and denote $\pi = w_1\cdots w_{n-2}$, then 
$\pi \in S_{n-1}$, hence $\pi=[b_1,\ldots ,b_{n-1},n]$.
If $w_{n-1}=1$ then $\sg \in S_{n-1}$ and we are done by induction. 
Hence assume $w_{n-1}\not =1$, so that $w_{n-1} =s_{n-1}\cdots s_k$ for some
$1\le k\le n-1$, and therefore 
$\sg =[b_1,\ldots ,b_{k-1},n,b_k,\ldots ,b_{n-1}]$.

\smallskip

{\bf Case 1:} $1\le k\le q$, in which case 
\[
\ell _q (w_{n-1})=n-q\quad\mbox{and}\quad 
\sg =[b_1,\ldots ,b_{k-1},n,b_k,\ldots ,b_q, \ldots , b_{n-1}].
\]
Then for $q\le i\le n-1$, 
\[
\#\{j<i+1\mid \sg(j)>\sg (i+1)\}=
\#\{j<i\mid b_j>b_i\}+1
\]
 (the ``$+1$'' comes from $n>b_i$).
It follows that 
$\,m_q(i+1,\sg )=m_q(i,\pi )+1$, since
\[
m_q(i+1,\sg )=\min \{i+1-q ;\#\{j<i+1\mid \sg (j)>\sg (i+1)\} \}=
\]
\[
=\min \{i+1-q ; \#\{j<i\mid b_j>b_i\}+1\}=
\]
\[
=\min \{i-q ; \#\{j<i\mid b_j>b_i\}\}+1=m_q(i,\pi )+1.
\]

\smallskip

Thus 
\[
inv_q(\sg)=\sum_{i=q+1}^nm_q(i,\sg)=\sum_{i=q}^{n-1}m_q(i+1,\sg)=
\sum_{i=q}^{n-1}m_q(i,\pi)+(n-q)=
\]
(by induction)
\[
=\ell_q(\pi)+n-q=\ell_q(\pi)+\ell_q(w_{n-1})=\ell_q(\sg).
\]

\smallskip

{\bf Case 2:} $q+1\le k$, ~~hence ~~$\ell_q(w_{n-1})=\ell_1(w_{n-1})=n-k$, 
$\sg =[b_1,\ldots ,b_q, \ldots ,b_{k-1},n,b_k, \ldots , b_{n-1}]$.
Here 
\begin{enumerate}
\item
$m_q(i,\sg)=m_q(i,\pi)\;$ if $\;q+1\le i \le k-1$,
\item
$m_q(k,\sg)=0$ $\;(i=k),\;$\\
\vskip 0.002 truecm
and, as in Case 1,
\item
$m_q(i+1,\sg)=m_q(i,\pi)+1\;$ if $\;k\le i\le n-1$.
\end{enumerate}
It follows that 
\[
inv_q(\sg)=\sum_{i=q+1}^nm_q(i,\sg)=\sum_{i=q+1}^{k-1}m_q(i,\pi)
+\sum_{i=k}^{n-1}m_q(i,\pi)+n-k=
\]
\[
\sum_{i=q}^{n-1}m_q(i,\pi)+(n-k)=
\]
(by induction)
\[
=\ell_q(\pi)+n-k=\ell_q(\pi)+\ell_q(w_{n-1})=\ell_q(\sg).
\]\qed

The following lemma was proved in~\cite{RR}.

\begin{lem}\label{RR-3.7}\cite[Lemma 3.7]{RR}
Let $w=s_{i_1}\cdots s_{i_p}$ be the  canonical presentation of
$w\in S_n$. Then the canonical presentation of $w^{-1}$ is obtained
from the presentation $w^{-1}=s_{i_p}\cdots s_{i_1}$ by commuting
moves only -- without any braid moves.
Similarly for $v,\,v^{-1}\in A_{n+1}$.
\end{lem}

\begin{pro}\label{od}
For every $\sg\in S_n$, 
\[
\ell_q(\sg ^{-1})=\ell_q(\sg ),\quad\mbox{hence also}\quad 
inv_q(\sg ^{-1})=inv_q(\sg ).
\]
\end{pro}
\noindent{\bf Proof.}  Lemma~\ref{RR-3.7} easily implies that
$\ell_q(\sg ^{-1})=\ell_q(\sg )$, while this, together with 
Proposition~\ref{cover01} imply the equality
$inv_q(\sg ^{-1})=inv_q(\sg )$.  \qed

\section{Basic Concepts II}\label{bas2}

A natural $q$-analogue of the $del$ statistics from~\cite{RR}
is introduced in this section. This allows us to introduce below a
(less intuitive) $q$-analogue of the descent statistics.

\subsection{The $del$ Statistics}\label{dels}

\noindent Recall the definitions of $Del$ and $del$ (of types $S$ and $A$)
from~\cite{RR}: Given a permutation $w$ in 
$S_n$, $Del_S(w)$ is the set of indices which
are {\it left to right minima} (l.t.r.min) in $w$, 
and $Del_A(w)$ is the set of indices which 
are {\it almost left to right minima} (a.l.t.r.min) in $w$. 
Let $s_i=(i,i+1), ~i=1,\ldots ,n-1$, denote the Coxeter generators
of $S_n$. The following classical fact is of fundamental importance
in this paper.

\smallskip

Let $R_j=\{ 1,s_j,s_js_{j-1},\ldots , s_js_{j-1}\cdots s_1\}$ and let
$w\in S_n$, then there exist unique elements $w_j\in R_j$,
 $1\le j\le n-1$, such that $w=w_1\cdots w_{n-1}$; this is the 
(unique) canonical presentation of $w$, see  Theorem 3.1 in~\cite{RR}. 

\smallskip

Similarly $a_i=s_1s_{i+1},~i=1,\ldots ,n-1$, 
are the corresponding generators for the alternating group $A_{n+1}$, and there
is a corresponding unique canonical presentation for the elements of $A_{n+1}$,
see Section 3 in~\cite{RR}.
The following was observed in~\cite{RR}: 
\begin{enumerate}
\item 
The number of times $s_1$ appears in the canonical presentation of
$w$ (i.e.~$del_S(w)$)
equals the number of l.t.r.min in  $w$ (hence
$\# Del_S(w)=del_S(w)$), see~\cite{RR} Proposition 7.2.
\item
The number of times $s_2$ appears in $w$ -- 
equals the number of a.l.t.r.min in  $w$. 
Moreover, if $w\in A_{n+1}$, that
number equals the number of times $a_1^{\pm 1}$ appears in 
the $A$-canonical presentation of $w$, which by 
definition is $del_A(w)$, and $del_A(w)=\# Del_A(w)$,
see~\cite{RR} Proposition 7.6.
\end{enumerate}

In this paper, `sub S' is replaced by `sub 1': 
$Del_S=Del_1$ and $del_S=del_1$,  etc.
Similarly (in $A_n$) `sub A' is replaced by `sub 2'.
We shall also encounter `sub $q$' for every positive integer $q$. 

\begin{df}\label{delq}
Let $\pi\in S_n$ and let $1\le q\le n-1$.
\begin{enumerate}
\item
Define $del_q(\pi)$ to be
the number of times $s_q$ appears in the canonical presentation of $\pi$.
\item
For $\,0\le k\le n-1$
define the $k$-th almost left to right minima
in a permutation $\pi\in S_n$ (denoted a$^k$.l.t.r.min)
as the set of indices
$$
Del_{k+1}(\pi):=
\{ i \mid k+2 \le i \le n, \#\{j<i \mid \pi (j) < \pi (i) \} \le k \}.
$$
Thus $Del_q(\pi)$ is the set of a$^{q-1}$.l.t.r.min in $\pi$.
\end{enumerate}
\end{df}
See Example~\ref{exa} below.

\medskip

{\bf Note} that if $i\le k+1 $ then, trivially, 
$\#\{j<i \mid \pi (j) < \pi (i) \} \le k$, however these indices are not
counted as a$^k$.l.t.r.min. Also note that a$^0$.l.t.r.min 
is simply  l.t.r.min.

\begin{pro}\label{altr2}
Let $w\in S_n$. Then for every nonnegative integer $k$,
the number of occurences of $s_{k+1}$ 
in the canonical presentation of $w$, $del_{k+1}(w)$,
equals the number of a$^k$.l.t.r.min in  $w$. Writing $k+1=q$ we have
\[
\# Del_q(w) = del_q(w).
\]
\end{pro}
\noindent{\bf Proof} (generalizes the proof of Proposition 7.6 in~\cite{RR}).\\
We first need the following two lemmas.
\begin{lem}\label{leper1} 
Let $1\le k+1\le n$, let $w\in S_n$
and let $\pi\in S_{k+1}$. Also let $i\le n$. Then $i$ is 
a$^k$.l.t.r.min of $w$ if and only if $i$ is 
a$^k$.l.t.r.min of $\pi w$. In particular, 
the number of  a$^k$.l.t.r.min of $w$ equals
the number of a$^k$.l.t.r.min of $\pi w$.
\end{lem}
\noindent{\bf Proof.} Denote $w=[b_1,\ldots , b_{n}]$ (namely $w(r)=b_r$),
and compare $w$ with $\pi w$:
$\pi $ permutes only the $b_r$'s in $\{1,\ldots ,k+1\}$. If 
$b_i\in\{1,\ldots ,k+1\}$, the total number of $b_j$'s smaller than 
$b_i$ is $\le k$; in particular such $i$ is a$^k$.l.t.r.min
in both $w$ and $\pi w$, provided
$i\ge k+2$. If on the other hand $b_i\not\in\{1,\ldots ,k+1\}$ then
$b_i$ is greater than all the elements in that subset; thus such
$i$ is a$^k$.l.t.r.min of $w$ if and only if $i$ is a$^k$.l.t.r.min of $\pi w$.
This implies the proof.
\qed
\begin{lem}\label{leper2}
Let $1\le k\le n-1$ and denote $s_{[k,n-1]}= s_ks_{k+1}\cdots s_{n-1}$.
Let $\sg\in S_{n-1}$ 
and write $\sg =[b_1,\ldots , b_{n-1},n]$. Then
$s_{[k,n-1]}\sg=[c_1,\ldots , c_{n-1},k]$, and the two tuples
$(b_1,\ldots , b_{n-1})$ and $(c_1,\ldots , c_{n-1})$ are order-isomorphic,
namely for all $i,j$, ~$b_i<b_j$ if and only if $c_i<c_j$. 
\end{lem}
\noindent{\bf Proof.}
Comparing $\sg $ with $s_{[k,n-1]}\sg $, we see that 
\begin{enumerate}
\item
the (position with) $n$ in $\sg $ is replaced in  $s_{[k,n-1]}\sg$  by $k$;
\item
each $j$ in $\sg $, $k\le j\le n-1$, is replaced by $j+1$ in $s_{[k,n-1]}\sg$;
\item
each $j$, $1\le j \le k-1$ is unchanged. 
\end{enumerate}
This implies the proof.   \qed

\medskip

\noindent{\bf The Proof of Proposition~\ref{altr2}} is by induction on $n$.
If $n\le k+1$, the number of a$^k$.l.t.r.min of any permutation in $S_n$ is
zero, and also $s_{k+1}\not\in S_n$, hence~\ref{altr2} holds
in that case.

Next assume~\ref{altr2} holds for $n-1$ and prove for $n$.  
Let $w=w_1\cdots w_{n-1}$ be the canonical presentation of $w\in S_n$ 
and denote $\sg =w_1\cdots w_{n-2}$, then $\sg \in S_{n-1}$. If $w_{n-1}=1$ 
then $w\in  S_{n-1}$ and the proof follows by induction. So let 
$w_{n-1}\not =1$, then we can write
$w_{n-1}=s_{n-1}s_{n-2}\cdots s_dv$, where 
$d\ge k+1$ and $v\in\{1,s_k,s_ks_{k-1},\ldots,s_ks_{k-1}\cdots s_1\}$ hence 
$v\in S_{k+1}$. If $d\ge k+2$ then necessarily $v=1$ and in that 
case the number
of times $s_{k+1}$ appears in $w$ and in $\sg$ is the same. If 
$d= k+1$, that number in $w$ is one more than in $\sg$. We show that the same
holds for the number of a$^k$.l.t.r.min for these two permutations
$\sg$ and $w$. 

By Lemma 3.4 of~\cite{RR}, it suffices to prove that statement
for the inverse permutations $w^{-1}$ and $\sg^{-1}$. Now,
$w ^{-1}=\pi s_{[d,n-1]}\sg^{-1}$, where $\pi = v^{-1}\in S_{k+1}$, hence by 
Lemma~\ref{leper1} it suffices to compare the number of a$^k$.l.t.r.min
in $\sg^{-1}$ with that in $s_{[d,n-1]}\sg^{-1}$. By Lemma~\ref{leper2}
$\sg^{-1}=[b_1,\ldots ,b_{n-1},n]$ and
$s_{[d,n-1]}\sg^{-1}=[c_1,\ldots ,c_{n-1},d]$ where
the $b$'s and the $c$'s are order isomorphic.

{\bf The case $d\ge k+2$}. Here the two last positions -- $n$ in $\sg^{-1}$
and $d$ in $s_{[d,n-1]}\sg^{-1}$ -- are not a$^k$.l.t.r.min, and the above
order isomorphism implies the proof in that case.

{\bf The case $d= k+1$}. By a similar argument, now the last position in
$s_{[d,n-1]}\sg^{-1}$ (which is $k+1$) is one additional a$^k$.l.t.r.min.\\
The proof now follows. 
\qed

\bigskip

\begin{pro}\label{del-inv}
For every positive integer $q$ and every permutation $\pi\in S_{n+q-1}$
$$
del_q(\pi)=del_q(\pi^{-1}).
$$
\end{pro}

\noindent{\bf Proof.}
This is a straightforward consequence of Lemma 3.7 of~\cite{RR}, 
which says the following: Let $\pi\in S_n$ and let
$\pi=s_{i_1}\cdots s_{i_r}$ be its canonical presentation. Then the canonical 
presentation of $\pi ^{-1}$ is obtained from the equation
$\pi ^{-1}=s_{i_r}\cdots s_{i_1}$ by commuting moves only, without
any braid moves. Thus, the number of times a particular $s_j$ appears in 
$\pi$ and in $\pi ^{-1}$ - is the same. This clearly implies the proof. \qed

\medskip
 
%
%
\begin{cor}\label{del-inv-cor}
For every positive integer $q$ and every permutation $\pi\in S_{n+q-1}$
the number of a$^{q-1}$.l.t.r.min in  $\pi$
equals the number of a$^{q-1}$.l.t.r.min in  $\pi^{-1}$.
\end{cor}

\noindent{\bf Proof.}
Combining Proposition~\ref{altr2}
with Proposition~\ref{del-inv}.
\qed

\begin{rem}\label{right1} 
Setting $q=k+1$ in Lemma~\ref{leper1}, deduce that 
for any two permutations $\sigma$  and $\eta$ in $S_{n+q-1}$,
if $\sigma$  and $\eta$ belong to the same right coset of $S_q$,
i.e.~$\eta\in S_q\sg$, then
$$
Del_q(\eta)=Del_q(\sigma)\qquad\mbox{and therefore}
\qquad del_q(\eta)=del_q(\sigma).
$$
The same is also true for the left cosets: 
Let $\eta\in \sg S_q$
then again
$$
Del_q(\eta)=Del_q(\sigma)\qquad\mbox{(and therefore}
\qquad del_q(\eta)=del_q(\sigma)).
$$
This easily follows from
Definition~\ref{delq}, since if $\sg =[b_1,\ldots ,b_q,\ldots ,b_n]$,
$\tau\in S_q$ and $\eta = \sg \tau$, then 
$\eta = [b_{\tau (1)},\ldots ,b_{\tau (q)},b_{q+1}\ldots ,b_n]$.

Let now $\sg$ and $\eta$ belong to the same left coset or right coset
of $S_q$, then by the same reasoning,
for any $q\le d$, $\,del_d(\eta)=del_d(\sigma)$ since 
$S_q\subseteq S_d$. Since 
\[
\ell_q (\eta)=\sum _{d=q}^{n-1}del_d(\eta),\qquad\mbox{and}\qquad
\ell_q (\sg)=\sum _{d=q}^{n-1}del _d(\sg),
\]
deduce that in that case $\ell_q (\eta)=\ell_q (\sg)$.
\end{rem}

\medskip

\subsection{The $q$-Descent Set}\label{qdes}

Recall that $i$ is a {\it descent} of $\pi$ if $\pi (i)>\pi (i+1)$, and
let $Des (\pi)$ denote the (`classical') descent-set of $\pi$.
The following definition seems to be the appropriate $q$-analogue for 
descents.

\begin{df}\label{id1} 
$i$ is a $q$-descent in $\pi\in S_{n+q-1}$ if $\ i\ge q$ and
at least one of the following two conditions holds:
$$
i\in Des(\pi) \leqno(1)
$$
$$
i+1 \hbox{ is an a$^{q-1}$.l.t.r.min in $\pi$}, \qquad \leqno(2)
$$
\end{df}

\medskip

\noindent {\bf Thus} 
$Des_q(\pi)=(Des(\pi)\cap \{q,q+1,\ldots ,n-1\})\cup (Del_q(\pi)-1)$,
hence for all $q$, ~$Del_q(\pi)-1\subseteq Des_q(\pi)$
where $Del_q(\pi)-1=\{i-1\mid i\in Del_q(\pi)\}$.

\noindent  Note that
when $q=1$, condition (2) says that $i+1$ is l.t.r.min, which implies that 
$i$ is a descent. Thus, a 1-descent is just a descent in the classical sense.

\begin{df}\label{id2} 
\begin{enumerate}
\item
The $q$-descent set of $\pi\in S_{n+q-1}$ is defined as
$$
Des_q(\pi):=\{i\ |\ i \hbox{ is a $q$-descent in } \pi\},   
$$
and 
\item
the $q$-descent number is defined as
$$
des_q(\pi):=\#Des_q(\pi). 
$$

\item
For $\pi\in S_{n+q-1}$ define
$$
maj_{q}(\pi):=\sum\limits_{i\in Des_q(\pi)} i
$$
and
$$
rmaj_{q,m}(\pi):=\sum\limits_{i\in Des_q(\pi)} (m-i),
$$
where $m=n+q-1$.
\end{enumerate}
\end{df}


\begin{exa}\label{exa}

Let $\sg=[7,8,6,5,2,9,4,1,3]$. 

When $q=2$, 
$Del_2(\sg )=\{3,4,5,7,8\}$ and 
$Des_2(\sg)=Del_2(\sg )-1=\{2,3,4,6,7\}$.

When $q=3$, $Del_3(\sg)=\{4,5,7,8,9\}$, hence 
$Des_3(\sg)=\{3,4,6,7\}\cup\{3,4,6,7,8\}=\{3,4,6,7,8\}$.

Also, $Des_4(\sg)=\{4,6,7,8\}$, etc.
\end{exa}

\section{Motivating Examples.}\label{ex}

When $q=1$, the corresponding statistics are classical. 
By definition,
for every $\pi\in S_n$
$$
\ell_1(\pi)=\ell_S(\pi),
$$
$$
\D_1(\pi)=\D_S(\pi), 
$$
and
$$
Del_1(\pi)=Del_S(\pi).
$$
It follows that
for every $\pi\in S_n$
$$
des_1(\pi)=des_S(\pi),
$$
$$
maj_1(\pi)=maj_S(\pi), \qquad ramj_{1,n}(\pi)=rmaj_{S_n}(\pi),
$$
and
$$
del_1(\pi)=del_S(\pi).
$$
The delent statistics, $del_S$, was introduced in~\cite{RR}.

The corresponding A-statistics were also studied in~\cite{RR}; these 
A-statistics correspond to the case $q=2$ -- and restricted to the 
alternating groups. This is the following proposition.

\begin{pro}\label{ex3}
For every even permutation $\pi\in S_{n+1}$
$$
\ell_2(\pi)=\ell_A(\pi), \leqno(1)
$$
$$
\D_2(\pi)=\D_A(\pi), \leqno(2)
$$
and
$$
Del_2(\pi)=Del_A(\pi). \leqno(3)
$$
\end{pro}

\noindent{\bf Proof.}
(1) follows from~\cite[Proposition 4.5]{RR}.

(2) follows from Lemma~\ref{geom}
in Appendix I (Section~\ref{appendixI}).

For (3) see~\cite[Proposition 7.5]{RR}.

\qed

\smallskip

An alternative and more conceptual proof is given below 
(see Remark~\ref{6.6-4.1}).

\medskip

\begin{cor}\label{ex4}
For every even permutation $\pi\in S_{n+1}$
$$
des_2(\pi)=des_A(\pi)
$$
$$
maj_2(\pi)=maj_A(\pi), \qquad ramj_{2,n}(\pi)=rmaj_{A_n}(\pi),
$$
and
$$
del_2(\pi)=del_A(\pi).
$$
\end{cor}

\section{The Double Cosets 
of $S_q\subseteq S_{n+q-1}$}\label{cosets}

Let $S_q$ be the subgroup of $S_{n+q-1}$ generated
by $\{s_1,\dots,s_{q-1}\}$. It is shown here that the previous
$q$-statistics are invariant on the double cosets of $S_q$
in $S_{n+q-1}$.


\begin{pro}\label{cover0}
For any two permutations $\pi$ and $\sigma$ in $S_{n+q-1}$,
if $\pi$ and $\sigma$ belong to the same double coset of $S_q$
(namely, $\pi\in S_q\sigma S_q$), then
\begin{enumerate}
\item
$$
Del_q(\pi)=Del_q(\sigma),\quad\mbox{hence}\quad del_q(\pi)=del_q;
$$
\item
$$
Des_q(\pi)=Des_q(\sigma),\quad\mbox{hence}\quad des_q(\pi)=des_q;
$$
\item
$$
inv_q(\pi)=inv_q(\sigma)=
\ell_q(\pi)=\ell_q(\sigma).
$$
\end{enumerate}
\end{pro} 
\noindent{\bf Proof.}
It suffices to prove that if there exists $\tau\in S_q$, such that
$\pi=\tau\sg$ or $\pi=\sg\tau$, then equalities 1--3 hold.

\smallskip

\noindent 1. Part 1 was proved in Remark~\ref{right1}.\\

\medskip

\noindent 2. Denote $\sg =[b_1,\ldots ,b_{n+q-1}]$ and
$\pi =[b'_1,\ldots ,b'_{n+q-1}]$.
Since
$
Des_q(\pi)=(Des(\pi)\cap \{q,q+1,\ldots ,n\})\cup (Del_q(\pi)-1),
$
and the same for $Des_q(\sg)$,
it suffices to prove the following:
Let $i\ge q$ and $i\in Des(\sg )$, then either $i\in Des(\pi)$ or
$i+1\in Del_q(\pi)$.

\smallskip
We prove first the case of the right cosets: $\pi=\tau\sg $.
It is given that $b_i > b_{i+1}$.\\
{\bf Case 1.} $b_i , b_{i+1}\not\in\{1,\ldots ,q\}$. Then $b_i=b'_i$
and $b_{i+1}=b'_{i+1}$ and we are done.\\
{\bf Case 2.} $b_i\not\in\{1,\ldots ,q\}$ and $b_{i+1}\in\{1,\ldots ,q\}$.
Then $b_i=b'_i>q$ while $b'_{i+1}\in\{1,\ldots ,q\}$ and we are done.\\
{\bf Case 3.} $b_i,b_{i+1}\in\{1,\ldots ,q\}$. Then at most $q-1$ $b_j$s 
in $\sg$
are left and smaller than $b_{i+1}$. Thus (by 1) 
$i+1\in Del_q(\sg )=Del_q(\pi )$. 

\medskip

We prove next the case of the left cosets: $\pi=\sg\tau $.\\
By the argument in Remark~\ref{right1}, the claim
holds if $i>q$. Therefore examine the case $i=q$. If $q\in Des(\pi)$,
then we are done. Recall that $b_q>b_{q+1}$ and assume 
$q\not\in Des(\pi)$ (i.e. $b_{\tau (q)}<b_{q+1}$). It follows that
\[
\#\{j<q+1\mid b_{\tau (j)}<b_{q+1}\}<q,
\]
hence $q+1\in Del_q(\pi)$, which completes the proof of part 2.

\noindent 3. This follows from Remark~\ref{right1} and from 
Proposition~\ref{cover01}, since
$inv_q(\pi)=\ell_q(\pi)$ and similarly for $\sg$. \qed

\bigskip

%
%
%
%
%
%
%
%

\section{The Covering Map $f_q$}\label{cover}

Motivated by Proposition~\ref{f2} below, we
introduce the map $f_q$ from $S_{n+q-1}$ onto $S_n$, 
which  sends all the elements
in the same double coset of $S_q$ to the same element in $S_n$.
The function $f_q$ is applied later to
``pull-back'' the equi-distribution results from
the (classical) case $q=1$ to the general $q$-case.

\bigskip

\begin{df}\label{covmap}
Let $\pi\in S_{n+q-1}$ and let $\pi=s_{i_1}\cdots s_{i_r}$ be its 
canonical presentation, then define $f_q:S_{n+q-1}\to S_n$ 
as follows:
\[
f_q(\pi)=f_q (s_{i_1})\cdots f_q (s_{i_r}),
\]
where $f_q (s_1)=\cdots =f_q (s_{q-1})=1$, and $f_q (s_j)=s_{j-q+1}$ if $j\ge q$.
\end{df}


\begin{rem}\label{cover01}
It is easy to verify that for any $q_1,q_2$, 
\[
f_{q_1}\circ f_{q_2}=f_{q_1+q_2-1}.
\]
Thus, for every natural $q$
$$
f_q=f_2^{q-1}.
$$
\end{rem}

%
\begin{pro}\label{dbl}
The map $f_q$ is invariant on the double cosets of $S_q$: Let 
$\sg\in S_{n+q-1}$ and $\pi\in S_q\sigma S_q$, then $f_q(\sg)=f_q(\pi)$.
\end{pro}
\noindent{\bf Proof.}
It suffices to prove that if $\sg\in S_{n+q-1}$ and
$\tau\in S_q$ then $f_q(\sg\tau)=
f_q(\tau\sg)=f_q(\sg)$.
By Remark~\ref{cover01}, it suffices to prove when $q=2$ and hence when 
$\tau=s_1$.
As usual, let $\sg=w_1\cdots w_{n}\in S_{n+1}$ be the canonical presentation of
$\sg$. By analyzing the two cases $w_1=1$ and $w_1=s_1$, it easily follows
that $f_2(s_1\sg)=f_2(\sg)$.

We now show that $f_2(\sg s_1)=f_2(\sg)$. 
The proof in that case follows
from the definition of $f_2$ and by induction on $n$, by analyzing the 
following cases:\\ 
$w_n=1$;\\ 
$w_n=s_ns_{n-1}\cdots s_k$ with $k\ge 3$;\\
$w_n=s_ns_{n-1}\cdots s_2$, and\\
$w_n=s_ns_{n-1}\cdots s_2s_1$. \\
We verify, for example, the case $k\ge 3$.
Denote $\pi=w_1\cdots w_{n-1}$, so  $\sg =\pi w_n$. Now 
$f_2(\sg s_1) =f_2(\pi s_1 \cdot w_n)=f_2(\pi s_1)f_2(w_n)=$
(by induction) $=f_2(\pi )f_2(w_n)=f_2(\sg)$.

The proof in the last two cases follows similarly, and from the fact that 
$f_2(s_ns_{n-1}\cdots s_2)=f_2(s_ns_{n-1}\cdots s_2s_1)
=s_{n-1}\cdots s_2s_1$.
\qed

\bigskip

Note that $f_q$ is not a group homomorphism. For example, let $q=2$,
$g=s_2$ and $h=s_1s_2$. Then $f_2(g)=f_2(h)=s_1$ so $f_2(g)f_2(h)=1$,
but $gh=s_1s_2s_1$, hence $f_2(gh)=s_1$. Nevertheless we do have the following 

\begin{pro}\label{hom}
For any permutation $\pi$, $f_q(\pi ^{-1})=(f_q(\pi))^{-1}$.
\end{pro}
\noindent{\bf Proof.} Again by Remark~\ref{cover01}, it suffices to prove for $q=2$.
The proof is based on Lemma~\ref{RR-3.7}. 

Denote $s_0:=1$, then note that if $s_is_j=s_js_i$ then also
$s_{i-1}s_{j-1}=s_{j-1}s_{i-1}$ (the converse is false, as $s_1s_2\not = 
s_2s_1$). 

Let $\pi=s_{i_1}\cdots s_{i_r} $ be the canonical presentation of $\pi$. By
commuting moves,
\[
\pi^{-1}=s_{i_r}\cdots s_{i_1}=\cdots =s_{p_1}\cdots s_{p_r} 
\]
where the right hand side is the canonical presentation of $\pi^{-1}$.
By definition, 
\[
f_2(\pi^{-1})=s_{p_1-1}\cdots s_{p_r-1}.
\]
Now by the same commuting moves
\[
s_{i_r-1}\cdots s_{i_1-1}=\cdots =s_{p_1-1}\cdots s_{p_r-1} 
\]
and the left hand side equals $(f_q(\pi))^{-1}$, which completes the proof. 
\qed

We also have
\begin{pro}\label{f2}
Recall from~\cite{RR} and Subsection~\ref{An} the map\\ 
$f: A_{n+1}\to S_n$. Then $f$ is the restriction $f=f_2|_{A_{n+1}}$
of $f_2$ to $A_{n+1}$.
\end{pro}

\noindent{\bf Proof.} Let $\pi\in A_{n+1}$, and let
\[
a_{i_1}^{\epsilon _1}\cdots a_{i_r}^{\epsilon _r}
\]
be its A--canonical presentation, where all $\epsilon _j=\pm 1$. By definition,
$f(\pi)=s_{i_1}\cdots s_{i_r}$.

Replace each $a_j$ in the above presentation by   $a_j=s_1s_{j+1}$ then,
by commuting moves `push' each $s_1$ as much as possible to the left. After 
some cancelations, an $s_1$ cannot move any more to the left if it is
already the left-most factor, or if it is preceded by an $s_2$ on its left.
It follows that 
\[
\pi=b s_{i_1+1}\cdots s_2s_1\cdots s_2s_1\cdots s_{i_r+1}\cdots
\]
where $b\in\{1,s_1\}$, and this is an S--canonical presentation. Then
$f_2(\pi)=s_{i_1}\cdots s_{i_r}$ and the proof follows.  \qed

\bigskip\bigskip

Restricting the maps  $f_q$ to 
$A_{n+q-1}$ -- we get more ``f-pairs'' (see~\cite{RR}, Section 5)
 with corresponding statistics,
equi-distributions and generating-functions-identities for the alternating
groups.

%
%

\bigskip

The main result here is
\begin{pro}\label{cover1}
For every $\pi\in S_{n+q-1}$
$$
Del_q(\pi)-q+1=Del_1(f_q(\pi)), \leqno(1)
$$
and in particular, $del_q(\pi)=del_1(f_q(\pi))$.
$$
Des_q(\pi)-q+1=Des_1(f_q(\pi)) \leqno(2)
$$
and in particular, $des_q(\pi)=des_1(f_q(\pi))$.
$$
inv_q(\pi)=inv_1(f_q(\pi))=\ell_q(\pi)=\ell_1 (f_q(\pi)). \leqno(3)
$$
Here $Del_q(\pi)-r=\{i-r\mid i\in Del_q(\pi)\}$
and similarly for $Des_q(\pi)-r$.

\end{pro}
The proof is given below.

\begin{rem}\label{DelM}
\begin{enumerate}
Recall that $R_j=\{1,s_j,s_js_{j-1},\ldots ,s_js_{j-1}\cdots s_1\}$.
\item
Let $w=w_1\cdots w_{n+q-2}$ where  all $w_j \in R_j$ be the 
canonical presentation of $w\in S_{n+q-1}$. Then 
$f_q(w)=f_q(w_1)\cdots f_q(w_{n+q-2})$ is the 
canonical presentation of $f_q(w)$. Note that 
$f_q(w_1)=\cdots =f_q(w_{q-1})=1$
\item
In addition, let also 
$w'=w'_1\cdots w'_{n+q-2}$, where also  $w'_j\in R_j$. It is obvious that 
$f_q(w)= f_q(w')$ if and only if $f_q(w_j)= f_q(w'_j)$ for all $j$.
\item
The definition of a$^{k}$.l.t.r.min in $\sg=[b_1,\ldots ,b_n]$ --
and therefore also the definition of the set
$Del_q(\sg)$ -- applies whenever 
the integers $b_1,\ldots ,b_n$ are distinct.
\item
Let $b_1,\ldots ,b_n$ and $c_1,\ldots ,c_n$ be two sets of distinct integers,
let $M$ be an integer satisfying $b_j,c_j<M$ for all $j$, let $1\le k\le n$
and denote
\[
\sg =[b_1,\ldots ,b_n],\qquad \sg^*=[b_1\ldots ,b_{k-1},M,b_k,\ldots ,b_n]
\]
and
\[
\eta =[c_1,\ldots ,c_n],\qquad \eta^*=[c_1\ldots ,c_{k-1},M,c_k,\ldots ,c_n].
\]
Then it is rather easy to verify that $Del_q(\sg)=Del_q(\eta)$ if and
only if $Del_q(\sg^*)=Del_q(\eta^*)$.
\end{enumerate}
\end{rem}
%

%
%
%
%
%
\begin{lem}\label{Delq2}
Let $w,w'\in S_{n+q-1}$ satisfy $f_q(w)= f_q(w')$, then 
\begin{enumerate}
\item
$Del_q(w)=Del_q(w')$.
\item
$Des_q(w)=Des_q(w')$.
\end{enumerate}
\end{lem}
\noindent{\bf Proof.} Since $f_1(w)=w$, we assume that $q\ge 2$.\\
\noindent 1.
By the definition of $f_q$ and by Remark~\ref{DelM}
it suffices to prove the following claim:

Let $w_j,w'_j\in R_j$ satisfy $f_q(w_j)= f_q(w'_j)$,
$q\le j\le n+q-2$, and let $w=w_q\cdots w_{n+q-2}$ and
$w'=w'_q\cdots w'_{n+q-2}$. Then $Del_q(w)=Del_q(w')$.

\smallskip

The proof is by induction on $n\ge 1$. If $n=1$, $w=w'=1$.\\
The induction step: 

Denote $m=n+q-1$, so $w=w_q\cdots w_{m-1}$ and 
$w'=w'_q\cdots w'_{m-1}$, then denote 
$\sg = w_q\cdots w_{m-2}$ and $\sg' = w'_q\cdots w'_{m-2}$. Since both
permutations are in $S_{m-1}\subseteq S_m$, we have
\[
\sg=[b_1,\ldots , b_{m-1},m]\qquad\mbox{and}\qquad
\sg'=[c_1,\ldots , c_{m-1},m].
\]
By induction, $Del_q(\sg)=Del_q(\sg')$. If $w_{m-1}=1$ then also $w'_{m-1}=1$
and we are done. 

Thus, assume both are $\not =1$.
Recall that $f_q(w_{m-1})=f_q(w'_{m-1})$ and
let $w_{m-1}=s_{m-1}\cdots s_k$ and $w'_{m-1}=s_{m-1}\cdots s_{k'}$.
If $k>q$, it follows that
$w_{m-1}=w'_{m-1}$ and we are done. So let $k,k'\le q$. 
By comparing both cases with the case $k=q$ we may assume that $k=q$
and $k'\le q$, hence $w'_{m-1}=w _{m-1}s_{q-1}\cdots s_{k'}$.

Compare  first $\sg w_{m-1}$ with $\sg' w_{m-1}$:
\[
\begin{array}{ll}
\sg w_{m-1}=[b_1,\ldots, b_{q-1},m,b_q,\ldots, b_{m-1}],\\
\sg' w_{m-1}=[c_1,\ldots, c_{q-1},m,c_q,\ldots, c_{m-1}],
\end{array}
\]
and by induction and Remark~\ref{DelM}.4,
$Del_q(\sg w_{m-1})=Del_q(\sg' w_{m-1})$.
Compare now $\sg' w_{m-1}$ with 
$\sg' w'_{m-1}=(\sg' w_{m-1})s_{q-1}\cdots s_{k'}$:
\[
\begin{array}{ll}
\sg' w_{m-1}=[c_1,\;\ldots\ldots\ldots\ldots\ldots , 
c_{q-1},m,c_q,\ldots\;, c_{m-1}]\qquad\mbox{and}\\
\sg' w'_{m-1}=[c_1,\ldots, c_{k'-1},m,c_{k'},\ldots,c_{q-1},
\ldots\ldots, c_{m-1}].
\end{array}
\]
A simple argument now shows that $\;q<i\;$ is  ~a$^{q-1}$.l.t.r.min 
 $Del_q(\sg' w_{m-1})=Del_q(\sg' w'_{m-1})$
and the proof of part 1 is complete.  

\medskip

\noindent 2. The proof is similar to that of Part 1.
Denote $m=n+q-1$, then
write $w = w_1\cdots w_{m-1}=\sg  w_{m-1}$ where
$\sg  = w_1\cdots w_{m-2}$, and similarly 
$w' = w'_1\cdots w'_{m-1}=\sg'  w'_{n-1}$. We assume that 
$f_q(w_j)=f_q(w'_j)$ for all $j$. Thus $f_q(\sg)=f_q(\sg ')$ and by 
induction, $Des_q(\sg)=Des_q(\sg')$. By an argument similar to that 
in the proof of part 1,
it follows that $Des_q(\sg w_{m-1})=Des_q(\sg' w_{m-1})$ and it
remains to show that $Des_q(\sg' w_{m-1})=Des_q(\sg' w'_{m-1})$. Again
as in the proof of part 1, we may assume that 
$w_{m-1}=s_{m-1}\cdots s_q$ and $w'_{m-1}=s_{m-1}\cdots s_t$ where $t<q$.
We prove the case $t=q-1$, the other cases being proved similarly.

Write $\sg' =[a_1,\ldots ,a_{m-1},m]$. Now 
$\sg' w'_{m-1}=\sg' w_{m-1}s_{q-1}$, hence
\[
\begin{array}{ll}
\sg' w_{m-1}=[a_1,\ldots ,a_{q-2},a_{q-1},m ,a_q,\ldots, a_{m-1}],\\
\sg' w'_{m-1}=[a_1,\ldots ,a_{q-2},m,a_{q-1},a_q,\ldots, a_{m-1}].
\end{array}
\]
Clearly, $q\in Des(\sg 'w_{m-1})$ (therefore $q\in Des_q(\sg 'w_{m-1})$), 
but it is possible that $q\not\in Des(\sg 'w'_{m-1})$. However, 
at most all the $q-1$ integers $a_1,\ldots , a_{q-1}$ are smaller
than $a_q$ (but $m > a_q$), hence $q+1\in Del_q(\sg'w'_{m-1})$, which implies
that $q\in Des_q(\sg'w'_{m-1})$. 

For all other indices $i\not = q$ it is easy to check that 
$i\in Des_q(\sg'w_{m-1})$ if and only if  $i\in Des_q(\sg'w'_{m-1})$, and
the proof is complete.  \qed

\bigskip

\noindent {\bf The Proof of Proposition~\ref{cover1}}. \\ 
1.
Let $\pi\in S_{n+q-1}$ and let
$\pi'$ denote the permutation obtained from $\pi$ 
by erasing -- in the canonical
presentation of $\pi$ -- all the appearances of the Coxeter generators 
$s_1,\ldots ,s_{q-1}$. Clearly, $f_q(\pi)=f_q(\pi')$, hence suffices to
prove that\\ 
(a)~ $Del_q(\pi)=Del_q(\pi')$, ~and \\
(b) ~$Del_q(\pi')-q+1=Del(f_q(\pi'))$,
i.e.  ~$Del_q(\pi')=Del(f_q(\pi'))+q-1$.\\

Let $\pi = w_1\cdots w_{q-1}w_q\cdots w_{m-1}$ ~($m=n+q-1$) be the canonical
presentation of $\pi$: ~$w_j\in R_j$. Denote $\tau=w_1\cdots w_{q-1}$
and $\sg=w_q\cdots w_{m-1}$, then both are given in their canonical 
presentations. Clearly, $f(\tau)=1$ and
$\pi'=\sg'=w'_q\cdots w'_{m-1}$, where 
for each $j$ ~$w'_j$ is obtained from $w_j$
by erasing all the appearances of $s_1,\ldots ,s_{q-1}$, and therefore
$f_q(w_j)=f_q(w'_j)$. 
By Lemma~\ref{Delq2}, 
$Del_q(\sg)=Del_q(\sg')=Del_q(\pi')$. Since $\pi=\tau\sg$ and 
$\tau\in S_q$, by Remark~\ref{right1} $Del_q(\pi)=Del_q(\sg)$ -- and (a)
is proved.

\smallskip

Part (b) follows from the following fact:\\ 
Let $\pi'=s_{i_1}\cdots s_{i_r}$ be the canonical presentation of 
the above $\pi'$
(therefore all $i_j \ge q$), then 
$f_q(\pi')=s_{i_1-q+1}\cdots s_{i_r-q+1}$. If
$f_q(\pi')=[a_1,\ldots , a_n]$, it follows that 
$\pi'=[1,\ldots, q-1,a_1+q-1,\ldots , a_n+q-1]$. If $2\le i$, it 
then follows that
$i$ is a l.t.r.min
of $f_q(\pi')$ if and only if $i+q-1$ is  ~~a$^{q-1}$.l.t.r.min
of $\pi'$, which proves (b).  \qed

%

\medskip

\noindent 2.
Recall that 
\[
Des_q(\pi)=(Des(\pi)\cap \{q,q+1,\ldots ,n\})\cup (Del_q(\pi)-1).
\]

{\it Special Case:} Assume $\pi$ does not involve any of $s_1,\ldots ,s_{q-1}$. 
As above, if
$f_q(\pi)=[a_1,\ldots , a_n]$ then 
$\pi=[1,\ldots, q-1,a_1+q-1,\ldots , a_n+q-1]$, hence
\[
Des(\pi)\cap \{q,q+1,\ldots ,n+q-1\}=Des(f_q(\pi))+q-1.
\]
By part 1
\[
Des_q(\pi)=\left ([Des(f_q(\pi))]\cup [Del(f_q(\pi))-1] \right )+q-1.
\]
Since for any $\sg\in S_n$ $Des(\sg)\supseteq Del(\sg)-1$, it follows that
the right hand side equals $Des(f_q(\pi))+q-1$, 
and this completes the proof of this case.

\smallskip

{\it The general case:} Let $\pi\in S_{n+q-1}$ be arbitrary. Let $\pi '$ be 
the permutation obtained  from $\pi$ by deleting all the appearances of 
$s_1,\ldots , s_{q-1}$ from its canonical presentation. Then 
$f_q(\pi)=f_q(\pi ')$ and the proof easily follows from the above special 
case and from Lemma~\ref{Delq2}(2).

\medskip

\noindent 3.
By Proposition~\ref{cover01}, $inv_q(\pi)=\ell _q(\pi)$. By the
definitions of $\ell _q$ and $ f_q $, ~$\ell _q(\pi)=\ell (f_q(\pi))$, 
and finally, $\ell (\sg)=inv (\sg)$ for any permutation $\sg$. \qed

\begin{rem}\label{6.6-4.1}
Proposition~\ref{ex3} now follows from Proposition~\ref{cover1},
combined with Proposition~\ref{f2} and \ref{RR-5.34}.
\end{rem}

\begin{lem}\label{cover2}
For every $\pi\in S_n$
$$
\# f_q^{-1}(\pi) =q!\cdot q^{del_1(\pi)}=(q-1)!\cdot q^{del_1(\pi)+1}.
$$
Moreover, let $g_q:A_{n+q-1}\to S_n$ be the restriction 
$g_q=f_q |_{A_{n+q-1}}$ of $f_q$ to $A_{n+q-1}$. Then 
\[
\# g_q^{-1}(\pi)=\frac{1}{2}\# f_q^{-1}(\pi).
\]
\end{lem}

\noindent{\bf Proof.} Denote $m=n+q-1$, so $f_q:S_m\to S_n$. Consider the canonical 
presentation of $\pi\in S_n$ and write it as $\pi=\pi^{(n-1)}\cdot v_{n-1}$, 
where $\pi^{(n-1)}\in S_{n-1}$ and $v_{n-1}\in R_{n-1}=
\{1,s_{n-1},s_{n-1}s_{n-2},\ldots ,s_{n-1}s_{n-2}\cdots s_1\}. $
Thus\\
 \[
\# f_q^{-1}(\pi)=\# f_q^{-1}(\pi^{(n-1)})\cdot
\# f_q^{-1}(v_{n-1})\ =
q!\cdot q^{del_1(\pi^{(n-1)})}\# f_q^{-1}(v_{n-1})
\]
 (by induction). If $del_1(v_{n-1})=0$ then $\# f_q^{-1}(v_{n-1})=1$.
If $del_1(v_{n-1})=1$ then $\# f_q^{-1}(v_{n-1})=q$, since in that case 
$v_{n-1}=s_{n-1}\cdots s_1$ and
\[
f_q^{-1}(v_{n-1})=\{w_{m-1},w_{m-1}s_{q-1},\ldots , w_{m-1}s_{q-1}\cdots s_1\},
\]
where $w_{m-1}=s_{m-1}s_{m-2}\cdots s_q$. The proof now follows. 

The argument for $g_q$ is similar. The factor $1/2$ comes from the fact
that $\#f_q^{-1}(1)=\#S_q$ while $\#g_q^{-1}(1)=\#A_q$. \qed

\medskip


Following~\cite{RR}, we introduce
\begin{df}\label{f-pairs1}
Let $m_1$ and $m_q$ be two statistics on the symmetric groups. We say that
$(m_1,m_q)$ is an  $f_q$--pair if for all $n$ and $\pi\in S_{n+q-1}$,
$m_q(\pi)=m_1(f_q(\pi))$.
\end{df}
As a corollary of Proposition~\ref{cover1} and Remark~\ref{iqu}, we have
\begin{cor}\label{f-pairs2}
The following are $f_q$--pairs:\\
$(inv_1,inv_q$), $(\ell_1,\ell_q)$, $(del_1,del_q)$,
$(des_1,des_q)$, and $(rmaj_{1,n},rmaj_{q,n+q-1})$.
\end{cor}
The same argument as in the proof of Proposition 5.6 in~\cite{RR}, 
together with Lemma~\ref{cover2}, now proves
\begin{pro}\label{f-pairs3}
Let $(m_1,m_q)$  be an  $f_q$--pair of statistics on the symmetric groups.
Then
\[
\sum_{\pi\in S_{n+q-1}}t_1^{m_q(\pi)}t_2^{del_q(\pi)}=
q!\sum_{\sg\in S_{n}}t_1^{m_1(\sg)}t_2^{del_1(\sg)}.
\]
Restricting $f_q$ to $A_{n+q-1}$ we obtain similarly, that
\[
\sum_{\pi\in A_{n+q-1}}t_1^{m_q(\pi)}t_2^{del_q(\pi)}=
\frac{1}{2}q!\sum_{\sg\in S_{n}}t_1^{m_1(\sg)}t_2^{del_1(\sg)}.
\]
\end{pro}

\begin{rem}\label{f-pairs4}
As in~\cite{RR}, Proposition~\ref{f-pairs3} allows us to lift 
equi-distribution theorems from $S_n$ to $S_{n+q-1}$, as well as to 
$A_{n+q-1}$. This is demonstrated in Theorem~\ref{qmac2}. We leave
the formulation and the proof of the corresponding
$A_{n+q-1}$ statement -- for the reader. 
\end{rem}



\medskip


\section{Dashed Patterns}\label{dash}

Dashed patterns in permutations were introduced in~\cite{BS}. For
example, the permutation $\sg$ contains the pattern $(1-32)$ if 
$\sg=[\ldots,a,\ldots,c,b,\ldots ]$ for some $a<b<c$; if no such
$a,b,c$ exist then $\sg$ is said to avoid  $(1-32)$. 
In~\cite{Cl} the author shows connections between the number of
permutations avoiding $(1-32)$ -- and various combinatorial objects, like 
the Bell and the Stirling numbers, as well as the number of
left-to-right-minima in permutations. 
In this and in the next sections we obtain the $q$-analogues
for these connections and results.

In Section~\ref{qdes} it was observed that, always, 
$Del_q(\pi )-1\subseteq Des_q(\pi )$. It is proved in Proposition~\ref{q-avoid}
that equality holds exactly for permutations avoiding a certain set of 
dashed-patterns.

\begin{df}\label{shem}
\begin{enumerate}
\item
Given $q$, denote by

\[
Pat(q)=\{(\sg_1-\sg_2-\cdots-\sg_q-(q+2)(q+1)) \mid \sg\in S_q\}
\]
the set with these $q!$ dashed patterns.\\
For example, $Pat(2)=\{(1-2-43),\;(2-1-43)\}$.
\item
Denote by $Avoid_q(m)$, $m=n+q-1$, the set of permutations in $S_{m}$
avoiding all the $q!$ patterns in $Pat(q)$, and let $h_q(m)$ denote
the number of the permutations in $S_{m}$  avoiding $Pat(q)$. Thus
$h_q(m)=\# Avoid_q(m)$ is the number of the permutations in 
$S_{n+q-1}$  avoiding $Pat(q)$.
Note that $h_q(m)=n!$ if $m\le q+1$. As usual, define $h_q(0)=1$.
%
\end{enumerate}
\end{df}
Connections between  $h_q(n)$ and the $q$-Bell and $q$-Stirling numbers 
are given in section~\ref{bell}.
\begin{rem}\label{sat}
A permutation $\pi\in S_{n+q-1}$ does satisfy one of the patterns in 
$Pat(q)$ if  and only if there exist a subsequence 
$$
1\le i_1<i_2<\dots<i_{q+1}< n+q-1,
$$
such that 
$
\pi(i_{q+1})>\pi(i_{q+1}+1)
$
and for every $1\le j\le q,\;$
$
\pi(i_j)<\pi(i_{q+1}+1).
$
In such a case, $i_{q+1}+1$ (namely, $\pi(i_{q+1}+1)$) is not
an a$^{q-1}$.l.t.r.min in $\pi$.
\end{rem}


\begin{pro}\label{q-avoid}
A permutation $\pi\in S_{n+q-1}$ avoids $Pat(q)$ exactly
when $Del_q(\pi)-1=Des_q(\pi)$:
$$
Avoid_q(n+q-1)=\{\pi\in S_{n+q-1}|\ Del_q(\pi)-1=Des_q(\pi)\}.
$$
In particular,
\[
h_q(n+q-1)=\# \{\pi\in S_{n+q-1}|\ Del_q(\pi)-1=Des_q(\pi)\} .
\]
\end{pro}
\noindent{\bf Proof.} 

1. Recall from Section~\ref{qdes} that, always, 
$Del_q(\pi )-1\subseteq Des_q(\pi )$. 
Let $\pi=[b_1,\ldots ,b_{n+q-1}]\in S_{n+q-1}$ satisfy 
$ Del_q(\pi)-1=Des_q(\pi)$, which implies that 
$Des(\pi)\cap\{q,\ldots ,n+q-1\}\subseteq Del_q(\pi)-1$,
and show that $\pi$ avoids $Pat(q)$. If not, by Remark~\ref{sat}
we obtain a descent in $\pi$  at $i_{q+1}$, while 
$i_{q+1}+1$ is not a$^{q-1}$.l.t.r.min in $\pi$; thus  
$i_{q+1}$ is in $ Des(\pi)\cap\{q,\ldots ,n+q-1\}$ 
but not in $Del_q(\pi)-1$, a contradiction.

\smallskip

2. Denote  $\pi=[b_1,\ldots ,b_{n+q-1}]$.
Assume now that $\pi\in Avoid_q(n)$, 
let $k\in Des(\pi)\cap\{q,\ldots ,n+q-1\}$ 
(so $b_k >b_{k+1}$) and show that $k+1\in Del_q(\pi)$, that is, $k+1$
(namely $b_{k+1}$) is a$^{q-1}$.l.t.r.min in $\pi$. If not, there 
exist $q$ (or more) $b_j$'s in $\pi$, smaller than and left 
of $b_{k+1}$ - hence also left of $b_{k}$. Together with $b_k >b_{k+1}$
this shows that $\pi\not\in Avoid_q(n+q-1)$, a contradiction.   \qed

\begin{cor}\label{qcor}
The covering map $f_q$ ~maps ~$Avoid_q(S_{n+q-1})$ ~to
 $Avoid _1(S_n)$:
\[
f_q: Avoid_q(S_{n+q-1})\to Avoid _1(S_n).
\]
Similarly,
\[
f_2: Avoid_q(S_{n+q-1})\to Avoid _{q-1}(S_{n+q-2}).
\]
\end{cor}
{\bf Proof.} This follows straightforward from Propositions~\ref{cover1}
and~\ref{q-avoid}. \qed

\section{$q$-Bell and $q$-Stirling Numbers}
\label{bell}
\subsection{The $q$-Bell Numbers}\label{qbell}

Recall that $S(n,k)$ are the Stirling numbers of the second kind, i.e.~the
numbers of $k$-partitions of the set $[n]=\{1,\ldots ,n\}$. Recall also that 
the Bell number $b(n)$ is the total number of the partitions of $[n]$:
$b(n)=\sum_kS(n,k)$.
\begin{df}\label{qbelln}
Define the $q$-Bell numbers $b_q(n)$ by 
\[
b_q(n)=\sum_k q^kS(n,k).
\]
\end{df}
\begin{rem}\label{wit}
Let $q\ge 1$
be an integer and consider partitions of $[n]$ 
into $k$ subsets,
where each subset is colored
by one of $q$ colors. The number of such $q$-colored
$k$-partitions is obviously $q^kS(n,k)$.
It follows that the total number of such 
$q$-colored partitions of  $[n]$ is the $n$-th $q$-Bell number $b_q(n)$.
\end{rem}

\noindent Proposition~\ref{qpro} below shows that
\[
\begin{array}{ll}
\#\{\sg\in S_{n+q-1}\mid Del_q(\sg)-1=Des_q(\sg)\quad\mbox{and}\quad
del_q(\sg)=k-1\} =~~~~~~~~~~~~~~~~~~\\
~~~~~~~~~~~~~~~~~~~~~~~~~~~~~~~~~~~~~~~~~~~~~~~~~~~~~~~~~~~~=
(q-1)!q^kS(n,k),
\end{array}
\]
and therefore
\[
(q-1)!b_q(n)=
\#\{\pi\in S_{n+q-1}|\ Del_q(\pi)-1=Des_q(\pi)\}.
\]
The $q$-Bell numbers are studied first.

\medskip

When $q=1$, by considering the subset in a $k$-partition of $[n]$
which contains $n$, one easily deduce the well-known recurrence relation 
\[
b_1(n)=\sum_k {n-1\choose k} b_1(n-k-1).
\]
In the general $q$ colors case,
apply the same argument, now taking into account that each 
subset -- and in particular the one containing $n$ -- can be colored by
$q$ colors. This proves:
\begin{lem}\label{rec}
 For each integer $q\ge1$ we have the following recurrence relation
\[
b_q(n)=q\sum_k {n-1\choose k} b_q(n-k-1).
\]
\end{lem}
\begin{rem}\label{remarkable}
\begin{enumerate}
\item 
Let $B_q(x)=\sum_{n=0}^{\infty}b_q(n)\frac{x^n}{n!}$ denote the 
exponential generating function  of $\{b_q(n)\}$.
As in page 42 in~\cite{W} , Lemma~\ref{rec} implies that 
$B'(x)=qe^xB_q(x)$. Together with 
$B(0)=1$ (since, by definition, $b_q(0)=1$), this implies that 
\[
B_q(x)=\exp (q e^x-q).
\]
\item
The classical formula 
\[
b_1(n)=\frac{1}{e}\sum_{r=0}^{\infty}\frac{r^n}{r!}
\]
~\cite{D}  (see also formula (1.6.10) in~\cite{W})
generalizes as follows: 
\[
b_q(n)=\frac{1}{e^q}\sum_{r=0}^{\infty}\frac{q^rr^n}{r!}.
\]
The proof follows, essentially unchanged, the argument on page 21 in~\cite{W}.
\end{enumerate}
\end{rem}

\subsection{Connections with Pattern-Avoiding Permutations}\label{qpat}

Recall that $Pat(q)=\{(\sg_1-\sg_2-\cdots-\sg_q-(q+2)(q+1)) \mid \sg\in S_q\}$
and that 
$h_q(n)$ denotes the number of the permutations in $S_n$  avoiding
all these $q!$ patterns in $Pat(q)$.
\begin{pro}\label{nu1}
The $q$-Bell numbers $b_q(n)$ and 
the numbers $\,h_q(n+q-1)$ of permutations in $S_{n+q-1}$ that avoid
$Pat(q)$, satisfy
\[
\,h_q(n+q-1)=(q-1)!\cdot b_q(n).
\]
By Proposition~\ref{q-avoid} this implies that
\[
(q-1)!b_q(n)=
\#\{\pi\in S_{n+q-1}|\ Del_q(\pi)-1=Des_q(\pi)\}.
\]
\end{pro}

The proof requires the following recurrence.
\begin{lem}\label{count}
If $n\ge q$ then
\[
h_q(n)=q\sum_{k=0}^{n-q}{n-q \choose k} h_q(n-k-1).
\]
\end{lem}
\noindent{\bf Proof.} The proof is by a rather standard argument.\\
Let $K\subseteq\{q+1,q+2,\ldots , n\}$ be a subset, with $|K|=k$, hence
$0\le k\le n-q$. Let $\kap$ be the word obtained by writing the numbers of $K$
in an increasing order. Note that there are ${n-q\choose k}$ such 
$K$'s -- hence ${n-q\choose k}$ such $\kap$'s. Let $1\le i\le q$ and let 
$\sg^{(i)}$ be a permutation of the set 
$\{1,\ldots ,i-1,i+1,\ldots ,n\}\setminus K$,
which avoids $Pat(q)$. 
By definition, since there are $n-1-k$ elements in that set, 
there are $h_q(n-k-1)$ such $\sg^{(i)}$'s. Now construct (the word)
$\eta^{(i)} = \sg^{(i)}i\kap$, then $\eta^{(i)}\in S_n$ and it avoids 
$Pat(q)$ since there is no descent in the part $i\kap$ of $\eta ^{(i)}$
(see Remark~\ref{sat}).
For each $1\le i\le q$, the number of 
$\eta^{(i)}$'s thus 
constructed is $\sum_{k=0}^{n-q}{n-q \choose k} h_q(n-k-1)$, 
hence 
\[
h_q(n)\le q\sum_{k=0}^{n-q}{n-q \choose k} h_q(n-k-1).
\]

Conversely, assume $\eta\in S_n$ avoids $Pat(q)$.
Among $1,\ldots , q$, let $i$ appear the rightmost in $\eta$ and write 
the word $\eta$ as  $\eta = \sg i \kap$, then none of $1,\ldots , q$
appears in $\kap$. The numbers in $\kap$ are increasing since otherwise, 
if there is a descent in $\kap$,
Remark~\ref{sat} would imply that $\eta$ does satisfy one of the dashed
patterns in $Pat(q)$, a contradiction. 
Since $\eta$ avoids $Pat(q)$, obviously the part $\sg$ of $\eta$ 
also avoids $Pat(q)$. It follows that $\eta$ is the above permutation
$\eta = \eta ^{(i)}$. This proves the reverse inequality -
and completes the proof.  \qed

\medskip

\noindent
{\bf The Proof of Proposition~\ref{nu1}}
now  follows by induction on $n\ge 0$. 

The case $n=0$ is clear. Assume $n\ge 1$, then by Lemma~\ref{count}
\[
 h_q(n+q-1)=q\sum_{k=0}^{n-1}{n-1 \choose k} h_q(n-1-k+q-1)=
\]
(by induction)
\[
= q\sum_{k=0}^{n-1}{n-1 \choose k}\cdot (q-1)!\cdot b_q(n-k-1)
=(q-1)!\cdot \left [q \sum_{k=0}^{n-1}{n-1 \choose k} b_q(n-k-1)\right ]
\]
\[
\mbox{(by Lemma~\ref{rec})}\qquad  =(q-1)!\cdot b_q(n).
\]  
This proves the first equation of the proposition.
Together with Definition~\ref{shem} and Proposition~\ref{q-avoid}, this 
implies that 
$h_q(n+q-1)=\#\{\pi\in S_{n+q-1}|\ Del_q(\pi)-1=Des_q(\pi)\}$,
hence
\[
(q-1)!b_q(n)=
\#\{\pi\in S_{n+q-1}|\ Del_q(\pi)-1=Des_q(\pi)\}.
\]
\qed

In the case $q=1$,
\[
b_1(n)=b(n)=\# Avoid_1(n)=
\# \{\sg\in S_n\mid Del_1(\sg ) -1 = Des_1(\sg )\}, 
\]
which appears in \cite{Cl}.

\bigskip

Let 
\[
H_q(x)=\sum_{n=0}^{\infty} h_q(n+q-1)\frac{x^n}{n!}
\]
be the exponential generating function  of the $h_q(n+q-1)$'s. 
By Remark~\ref{remarkable}(1) and 
Proposition~\ref{nu1} we have

\begin{cor}\label{}
\[
H_q(x)=(q-1)!\cdot \exp (q e^x -q ).
\]
\end{cor}

%

\subsection{Stirling Numbers of the Second Kind}\label{avq}

The following refinement of 
the second equation of Proposition~\ref{nu1}
is proved in this subsection.

\begin{pro}\label{qpro}
\[
\begin{array}{ll}
\#\{\sg\in S_{n+q-1}\mid Del_q(\sg)-1=Des_q(\sg)\quad\mbox{and}\quad
del_q(\sg)=k-1\} =~~~~~~~~~~~~~~~~~~\\
~~~~~~~~~~~~~~~~~~~~~~~~~~~~~~~~~~~~~~~~~~~~~~~~~~~~~~~~~~~~=(q-1)!q^kS(n,k).
\end{array}
\]
Deduce that\\
$$
\sum\limits_{\{\pi\in S_n|\ Del_1(\pi)-1=Des_1(\pi)\}} q^{del_1(\pi)}=
\frac{1}{q}\cdot b_q(n),
$$
and more generally,
\[
\begin{array}{ll}
\sum\limits_{\{\sg\in S_{n+q-1}|\ Del_q(\sg)-1=Des_q(\sg)\}} 
q^{del_q(\sg)}=~~~~~~~~~~~~~~~~~~~~~~~~~~~~~~~~~~\\
~~~~~~~~~~~~~~~~~~~~~~~~~~~~~~~~~=\frac{(q-1)!}{q}\cdot \sum_k q^{2k}S(n,k)=
\frac{(q-1)!}{q}\cdot b_{q^2}(n).
\end{array}
\]
\end{pro}
%
%

\noindent{\bf Proof.} We first prove the case $q=1$ namely, that
\[
\#\{\sg\in S_n \mid Del_1(\sg ) -1 = Des_1(\sg )\quad
\mbox{and}\quad del_1(\sg )=k-1\}=S(n,k).
\]
Recall that $S(n,k)$ is the number of partitions of $[n]$
into $k$ non-empty subsets. Given such a partition
$D_1\cup\cdots\cup D_k=[n]$, assume w.l.o.g.~that the numbers in each $D_i$
are increasing : $D_i$ is
$\{d_{i,1} <d_{i,2}<\cdots\}$, and also, the minimal elements 
$d_{1,1},d_{2,1},\cdots $ are decreasing: 
$d_{1,1}>d_{2,1}>\cdots >d_{k,1}$. Corresponding to that partition we 
construct the permutation $\sg =[D_1,D_2,\ldots ]$, namely
$\sg = [d_{1,1},d_{1,2},\ldots , d_{2,1},d_{2,2},
\ldots ,d_{k,1},d_{k,2}\ldots]$.

Now $Del_1(\sg )$, the l.t.r.min of $\sg$, are exactly at the $(k-1)$ 
positions of $d_{2,1},d_{3,1},\ldots ,d_{k,1}$, and obviously the 
descents occur at  $Del_1(\sg )-1$. This establishes an injection of
the set of the $k$ partitions of $[n]$ into the above set, which implies
that
\[
\mbox{card}\,\{\sg\in S_n \mid Del_1(\sg ) -1 = Des_1(\sg )\quad
\mbox{and}\quad del_1(\sg )=k-1\}\ge S(n,k).
\]
Since the sum on all $k$'s of both sides equals $b(n)$, this implies the 
case $q=1$.

\medskip

The general $q$ case follows from Proposition~\ref{cover1}, 
and from Lemma~\ref{cover2}: 

\noindent Let $\pi\in S_n$. By Proposition~\ref{cover1}, 
\[
Del_1(\pi)-1=Des_1(\pi)\quad\mbox{ if and only if}\quad 
Del_q(f_q^{-1}(\pi))-1=Des_q(f_q^{-1}(\pi)),
\]
and also, 
\[
del_1(\pi)=k-1\quad\mbox{ if and only if}\quad
del_q(f_q^{-1}(\pi))=k-1.
\]
Denote $D_q(n,k)=
\{\sg\in S_{n+q-1}\mid Del_q(\sg)-1=Des_q(\sg)\;\&\;
del_q(\sg)=k-1\}$, so that 
$D_1(n,k)=\{\pi\in S_{n}\mid Del_q(\pi)-1=Des_1(\pi)\;\&\;
del_1(\pi)=k-1\}$.  It follows that 
\[
D_q(n,k)=\bigcup _{\pi\in D_1(n,k)}f_q^{-1}(\pi),
\]
a disjoint union. By Lemma~\ref{cover2}, $\# f_q^{-1}(\pi)
=(q-1)!\cdot q^k$ for all $\pi\in D_1(n,k)$, and the proof now follows
easily from the case $q=1$. \qed


\subsection{Stirling Numbers of the First Kind}\label{firstk}

Let $c(n,k)$ 
be the signless Stirling numbers of the first kind.

\begin{pro}\label{qc1} 
$$
c(n,k)=\#\{\pi\in S_n|\ (del_S(\pi))=del_1(\pi)=k-1\} ,
$$
namely, $c(n,k)$ equals the number of permutations in $S_n$ with $k-1\;$
l.t.r.min.
\end{pro}
For the proof, see Proposition 5.8 in~\cite{RR}.

\medskip

The following is a $q$-analogue of Proposition \ref{qc1}.

\begin{pro}\label{qc3}
$$
\#\{\pi\in S_{n+q-1}|\ del_q(\pi)=k-1\}=c_q(n,k),
$$
where $c_q(n,k)=q^k\ (q-1)!\ c(n,k)$.
\end{pro}
\noindent{\bf Proof.} The proof is essentially identical to the proof of 
Proposition~\ref{qpro}, with the set $D_q(n,k)$ being replaced here by 
the set $H_q(n,k)=\{\pi\in S_{n+q-1}|\ del_q(\pi)=k-1\}$. Then 
$H_1(n,k)=\{\pi\in S_{n}|\ del_1(\pi)=k-1\}$, and by
Proposition 5.8 in~\cite{RR}, ~~$\# H_1(n,k) =c(n,k)$, the {\it signless
Stirling number of the first kind}. The proof now follows.  \qed

\medskip

%


\section{Equidistribution}\label{equid}

\subsection{MacMahon Type Theorems for $q$-Statistics}
\label{q-MM}



Recall the definition of $rmaj_{q,n+q-1}$ from Definition~\ref{id2}.

\begin{rem}\label{iqu}
Note that for $\pi\in S_{n+q-1}$, 
\[
rmaj_{q,n+q-1}(\pi)=rmaj_{1,n}(f_q(\pi))=rmaj_{S_n}(f_q(\pi)).
\]
This follows since by Proposition~\ref{cover1}(2), $i\in Des_q(\pi) $
if and only if $i-q+1\in Des_1(f_q(\pi))$.
\end{rem}

The following is a $q$-analogue of MacMahon's equi-distribution theorem.

\begin{thm}\label{qmac}
For every positive integers $n$ and $q$
$$
\begin{array}{ll}
\sum\limits_{\pi\in S_{n+q-1}} t^{rmaj_{q,n+q-1}(\pi)}
=
\sum\limits_{\pi\in S_{n+q-1}} t^{inv_q(\pi)}=~~~~~~~~~\\
~~~~~~~~~~~~~~~~~~~~~~=q!(1+tq)(1+t+t^2q)\cdots 
(1+t+\ldots +t^{n-2}+t^{n-1}q).
\end{array}
$$
\end{thm}
This theorem is obtained from the next one by substituting $t_2=1$.
\begin{thm}\label{qmac2}
For every positive integers $n$ and $q$
$$
\begin{array}{ll}
\sum\limits_{\pi\in S_{n+q-1}} t_1^{rmaj_{q,n+q-1}(\pi)}t_2^{del_q(\pi)}
=
\sum\limits_{\pi\in S_{n+q-1}} t_1^{inv_q(\pi)}t_2^{del_q(\pi)}=~~~~~~~~~\\
~~~~~~~~~~~=q!(1+t_1t_2q)(1+t_1+t_1^2t_2q)\cdots 
(1+t_1+\ldots +t_1^{n-2}+t_1^{n-1}t_2q).
\end{array}
$$
\end{thm}
\noindent{\bf Proof.}
By Proposition \ref{cover1} and Remark~\ref{iqu}, $(rmaj_{S_n}, rmaj_{q,n+q-1})$ and $(inv, inv_q)$
are $f_q$-pairs.
The proof now follows from Proposition~\ref{f-pairs3} and Theorem~\ref{6.1}. 

\qed

The following is a $q$-analogue of Foata-Sch\"utzenberger's equi-distribution theorem \cite[Theorem 1]{FS}.

\begin{thm}\label{fs-q}
For every positive integers $n$ and $q$ and
every subset $B\subseteq[q,n+q-1]$
$$
\sum\limits_{\{\pi\in S_{n+q-1}|\ Des_q(\pi^{-1})=B\}}
t^{inv_q(\pi)}=
\sum\limits_{\{\pi\in S_{n+q-1}|\ Des_q(\pi^{-1})=B\}}
t^{rmaj_{q,n+q-1}(\pi)}.
$$
\end{thm}

This theorem is obtained from the next one by substituting $B_2=[q,n+q-1]$.

\begin{thm}\label{fs-q2}
For every positive integers $n$ and $q$ and
every subsets $B_1\subseteq[q,n+q-1]$ and $B_2\subseteq [q,n+q-1]$
$$
\sum\limits_{\{\pi\in S_{n+q-1}|\ Des_q(\pi^{-1})=B_1,\ Del_q(\pi^{-1})=B_2\}}
t^{inv_q(\pi)}=
$$
$$
\sum\limits_{\{\pi\in S_{n+q-1}|\ Des_q(\pi^{-1})=B_1,\ Del_q(\pi^{-1})=B_2\}}
t^{rmaj_{q,n+q-1}(\pi)}.
$$
\end{thm}

\noindent{\bf Proof.}
By Proposition~\ref{cover1} and Remark~\ref{iqu}, 
it suffices to prove that for every subsets $B_1\subseteq [n-1]$
and $B_2\subseteq [n-1]$
$$
\sum\limits_{\{\pi\in S_{n+q-1}|\ Des_1(f_q(\pi^{-1}))=B_1,\ 
Del_1(f_q(\pi^{-1}))=B_2\}}
t^{inv_1(f_q(\pi))}=
$$
$$
=\sum\limits_{\{\pi\in S_{n+q-1}|\ Des_1(f_q(\pi^{-1}))=B_1,\ 
Del_1(f_q(\pi^{-1}))=B_2\}}
t^{rmaj_{1,n}(f_q(\pi))}.
$$
By Proposition~\ref{hom} $f_q(\pi^{-1})=f_q(\pi)^{-1}$.
Thus, denoting $\sg=f_q(\pi)$, it suffices to prove that
$$
\sum\limits_{\{\sg\in S_{n}|\ Des_1(\sg^{-1})=B_1,\ 
Del_1(\sg^{-1})=B_2\}} 
\#f_q^{-1}(\sg)\cdot
t^{inv_1(\sg)}=
$$
$$
\sum\limits_{\{\sg\in S_{n}|\ Des_1(\sg^{-1})=B_1,\ 
Del_1(\sg^{-1})=B_2\}} 
\#f_q^{-1}(\sg)\cdot
t^{rmaj_{1,n}(\sg)}.
$$
By Propositions~\ref{altr2} and~\ref{del-inv}, for every $\sg\in S_n$
with $Del_1(\sg^{-1})=B_2$ , $del_1(\sg)=\#B_2$.
Thus, by Lemma~\ref{cover2}, $\#f_q^{-1}(\sg)=(q-1)!\cdot q^{\#B_2+1}$
for all permutations in the sums. Hence, the theorem is reduced to
$$
(q-1)!\cdot q^{\#B_2+1}\cdot
\sum\limits_{\{\sg\in S_{n}|\ Des_1(\sg^{-1})=B_1,\ 
Del_1(\sg^{-1})=B_2\}} 
t^{inv_1(\sg)}=
$$
$$
(q-1)!\cdot q^{\#B_2+1}\cdot
\sum\limits_{\{\sg\in S_{n}|\ Des_1(\sg^{-1})=B_1,\ 
Del_1(\sg^{-1})=B_2\}} 
t^{rmaj_{1,n}(\sg)}.
$$
Theorem~\ref{9.1} completes the proof. 

\qed

\medskip

Theorem~\ref{fs-q} implies $q$-analogues of two classical
identities, due to~\cite{Ros, FS}.

\begin{cor}\label{q-rosel}
For every positive integers $n$ and $q$
$$
\leqno(1)
\sum\limits_{\pi\in S_{n+q-1}} t_1^{inv_q(\pi)} t_2^{des_q(\pi^{-1})}=
\sum\limits_{\pi\in S_{n+q-1}} t_1^{rmaj_{q,n+q-1}(\pi)} 
t_2^{des_q(\pi^{-1})},
$$
and
$$
\leqno(2)
\sum\limits_{\pi\in S_{n+q-1}} t_1^{inv_q(\pi)} 
t_2^{rmaj_{q,n+q-1}(\pi^{-1})}=
\sum\limits_{\pi\in S_{n+q-1}} t_1^{rmaj_{q,n+q-1}(\pi)} 
t_2^{rmaj_{q,n+q-1}(\pi^{-1})}.
$$
\end{cor}

\subsection{Equi-distribution on $Avoid_q(n)$}
\label{equid-avoid}

The main theorem on equi-distribution on permutations avoiding patterns
is the following.

\begin{thm}\label{q6}
For every positive integers $n$ and $q$ and
every subset $B\subseteq [q,\dots,n+q-2]$
$$
\sum\limits_{\{\pi ^{-1}\in Avoid_q(n+q-1)|\ Des_q(\pi ^{-1})=B\}} t^{rmaj_{q,n+q-1}(\pi)}
=
\sum\limits_{\{\pi ^{-1}\in Avoid_q(n+q-1)|\ Des_q(\pi ^{-1})=B\}} t^{inv_q(\pi)}
$$
\end{thm}

\bigskip

\noindent{\bf Proof.}
Substituting $B_1=B_2-1=B$ in Theorem~\ref{fs-q2} we obtain,
for every subsets $B\subseteq[q,n+q-1]$ 
$$
\sum\limits_{\{\pi\in S_{n+q-1}|\ Des_q(\pi^{-1})=Del_q(\pi^{-1})-1=B\}}
t^{inv_q(\pi)}=
$$
$$
\sum\limits_{\{\pi\in S_{n+q-1}|\ Des_q(\pi^{-1})=Del_q(\pi^{-1})-1=B\}}
t^{rmaj_{q,n+q-1}(\pi)}.
$$
By Proposition~\ref{q-avoid}
$$
\{\pi\in S_{n+q-1}|\ Des_q(\pi^{-1})=Del_q(\pi^{-1})-1=B\}
$$
$$
=\{\pi ^{-1}\in Avoid_q(n+q-1)|\ Des_q(\pi ^{-1})=B\}.
$$
This completes the proof.

\qed

%


%

Theorem~\ref{q6} implies

\begin{cor}\label{qs5}
For every positive integers $n$
and $q$
$$
\sum\limits_{\pi ^{-1}\in Avoid_q(n+q-1)} t_1^{rmaj_{q,n+q-1}(\pi)}t_2^{des_q(\pi)}
=
\sum\limits_{\pi ^{-1}\in Avoid_q(n+q-1)} t_1^{inv_q(\pi)}t_2^{des_q(\pi)}
$$
\end{cor}

%

\medskip

The following is an
extension of MacMahon's theorem to permutations avoiding patterns.

\begin{thm}\label{q5}
For every positive integers $n$ and $q$
$$
\sum\limits_{\pi ^{-1}\in Avoid_q(n+q-1)} t^{rmaj_{q,n+q-1}(\pi)}
=
\sum\limits_{\pi ^{-1}\in Avoid_q(n+q-1)} t^{inv_q(\pi)}
$$
\end{thm}

\noindent{\bf Proof.}
Substitute $t_2=1$ in Corollary~\ref{qs5}.
\qed

%






\section{Appendix. $Des_2=Des_A$: The Proof}
\label{appendixI}

\begin{lem}\label{geom} 
 Let $w=[b_1,\ldots ,b_{n+1}]\in A_{n+1}$. 
Let $1\le i \le n-1$, then
$i\in Des_A(w)$ if and only if one of the following two conditions hold.
\begin{enumerate}
\item
 $b_{i+1}>b_{i+2}$,\\ 
or
\item
 $b_{i+1}<b_{i+2}$ 
and $b_1,b_2,\ldots , b_i >b_{i+2}$. 
\end{enumerate}
In particular, $1\in Des_A(w)$ if and only if 
$b_1 > b_3$  (and/) or $b_2 > b_3$.
\end{lem}
\noindent{\bf Proof.} 
The basic tool is the formula
\[
\ell_A(w)=\ell_S(w)-del_S(w).
\]
Assume first that $2\le i \le n-1$, then 
$v=wa_i=[b_2,b_1,\ldots ,b_{i+2},b_{i+1},\ldots ]$.
Now compare $\ell_S(w)$ with $\ell_S(v)$, and $del_S(w)$ with $del_S(v)$, 
then apply
the above formula, and the proof follows. Here are the details.

\smallskip

{\bf The case} $2\le i \le n-1$ and $b_{i+1}>b_{i+2}$.\\
If $b_1<b_2$ then $\ell_S(w)=\ell_S(v)$. Now, $del(\sg)$ is the number
of l.t.r.min in $\sg$. Interchanging $b_1<b_2$ in $w$ adds one such
l.t.r.min, while interchanging $b_{i+1}>b_{i+2}$ reduces that ($del_S$)
number by one, or leaves it unchanged. In particular, 
$del_S(w) \le del_S(v)$. It follows that\\ 
$\ell_A(w)=\ell_S(w)-del_S(w)\ge \ell_S(v)-del_S(v)=\ell_A(v)$,
i.e.~$\ell_A(wa_i)\le\ell_A(w)$,
hence $i\in Des_A(w)$. 

\smallskip

Similarly for the other cases. If $b_1>b_2$ (and $b_{i+1}>b_{i+2}$), verify
that $\ell_S(w)=\ell_S(v)+2$, while $del_S(w)\le del_S(v)+2$, and again this
implies that $i\in Des_A(w)$.
This completes the proof of 2.a.

\smallskip

{\bf The case} $2\le i \le n-1$ and $b_{i+1}<b_{i+2}$.\\
First, assume $b_1<b_2$, then $\ell_S(v) =\ell_S(w)+2$. If 
$b_1,b_2,\ldots , b_i >b_{i+2}$ then also $del_S(v) =del_S(w)+2$,
hence $\ell_A(wa_i)=\ell_A(v)=\ell_A(w)$, and $i\in Des_A(w)$.
If  $b_j<b_{i+2}$ for some $1\le j\le i$ then $del_S(v) =del_S(w)+1$
and it follows that $i\not\in Des_A(w)$.

If $b_1>b_2$ then $\ell_S(v) =\ell_S(w)$. Assuming that 
$b_1,b_2,\ldots , b_i >b_{i+2}$, deduce that also $del_S(v) =del_S(w)$,
hence $i\in Des_A(w)$.
If  $b_j<b_{i+2}$ for some $1\le j\le i$ then $del_S(v) =del_S(w)-1$,
so  $\ell_A(wa_i)=\ell_A(v)=\ell_A(w)-1$ and $i\not\in Des_A(w)$.

\medskip
Finally assume that $i=1$, then 
$v=wa_1=ws_1s_2=[b_2,b_3,b_1,b_4,b_5\ldots ]$.
Obviously, $\ell_S(w)-\ell_S(v)$ depends only on the order relations among
$b_1, b_2,b_3$, and similarly for  $del_S(w)-del_S(v)$. We can therefore 
assume that $\{b_1, b_2,b_3\}=\{1, 2,3\}$, then check the 3!=6 possible
cases of $w=[b_1, b_2,b_3,\ldots ]$. 
For example, assume 
$w=[1,3,2,\ldots ]$, then $wa_1=[3,2,1,\ldots ]=v$, so $\ell_S(v) =\ell_S(w)+2$
while $del_S(v) =del_S(w)+2$, hence $1\in Des_A(w)$.\\
Similarly for the remaining  five cases. \qed

%
%

%

\end{document}